\documentclass[journal,draftcls,10pt,onecolumn]{IEEEtran}

\usepackage{amsfonts}
\usepackage{amsmath}
\usepackage{graphicx}

\newcommand{\R}{\mathbb{R}}
\renewcommand{\S}{\mathbb{S}}
\newcommand{\N}{\mathbb{N}}
\newcommand{\bmat}[1]{\begin{bmatrix}#1\end{bmatrix}}
\newcommand{\norm}[1]{\left\lVert{#1}\right\rVert}

\newcommand{\ip}[2]{\left\langle #1, #2 \right\rangle}
\newcommand{\half}{\frac{1}{2}}

\newtheorem{thm}{Theorem}

\newtheorem{lem}[thm]{Lemma}

\let\bbl\Bigl

\let\bbr\Bigr

\begin{document}

\title{\LARGE \bf
SOS Methods for Multi-Delay Systems:\\ A Dual Form of Lyapanov-Krasovskii Functional
}

\author{Matthew~M.~Peet,~\IEEEmembership{Member,~IEEE,}%
\thanks{M. Peet is with the School for the Engineering of Matter, Transport and Energy, Arizona State University, Tempe, AZ, 85298 USA. e-mail: (\tt \small mpeet@asu.edu}}

%

\maketitle

\begin{abstract}
We present a dual form of Lyapunov-Krasovskii functional which allows the problem of controller synthesis of multi-delay systems to be formulated and solved in a convex manner. First, we give a general form of dual stability condition formulated in terms of Lyapunov operators which are positive, self-adjoint and preserve the structure of the state-space. Second, we provide a class of such operators and express the stability conditions as positivity and negativity of quadratic Lyapunov-Krasovskii functional forms. Next, we adapt the SOS methodology to express positivity and negativity of these forms as LMIs, describing a new set of polynomial manipulation tools designed for this purpose. Finally, we apply the resulting LMIs to a battery of numerical examples and demonstrate that the stability conditions are not conservative. The results of this paper are significant in that they open the way for dynamic output $H_\infty$ optimal control of systems with multiple time-delays.
\end{abstract} 

\section{Introduction}
Systems with delay have been well-studied for some time~\cite{niculescu_book,gu_2003,richard_2003}. Recently, there have been many results on the use of optimization and semidefinite programming for stability of linear and nonlinear time-delay systems. Although the computational question of stability of a linear state-delayed system is believed to be NP-hard, several techniques have been developed which use LMI methods~\cite{boyd_book} to construct sequences of polynomial-time algorithms which provide sufficient stability conditions and appear to converge to necessity as the complexity of the algorithms increase. Examples of such sequential algorithms include the piecewise-linear approach~\cite{gu_2003}, the delay-partitioning approach~\cite{gouaisbaut_2009}, the Wirtinger-based method of~\cite{seuret_2013} and the SOS approach~\cite{peet_2009SICON}. In addition, there are also frequency-domain approaches such as~\cite{michiels_2005,sipahi_2005}. These algorithms are sufficiently reliable so that for the purposes of this paper, we may consider the problem of stability analysis of linear discrete-delay systems to be solved.

The purpose of this paper is to explore methods by which the success in stability analysis of time-delay systems may be used to attack what may be considered the relatively underdeveloped field of robust and optimal controller synthesis. Although there have been a number of results on controller synthesis for time-delay systems~\cite{luo_book}, none of these results has been able to resolve the fundamental bilinearity of the synthesis problem. That is, controller synthesis is not convex in the combined Lyapunov operator $\mathcal{P}$ and feedback operator $\mathcal{K}$. Without convexity, it is difficult to construct provably stabilizing controllers without significant conservatism, much less address the problems of robust and quadratic stability. Some papers use iterative methods to alternately optimize the Lyapunov operator and controller as in~\cite{moon_2001} or~\cite{fridman_2002} (via a ``tuning parameter''). However, this iterative approach is not guaranteed to converge. Meanwhile, approaches based on frequency-domain methods, discrete approximation, or Smith predictors result in controllers which are not provably stable or are sensitive to variations in system parameters or in delay. Finally, we mention that delays often occur in both state and input and to date most methods do not provide a unifying formulation of the controller synthesis problem with both state and input delay.

In this paper, we create a unified inequality-based framework for robust and optimal control of systems with multiple delays. The model for our approach is the LMI framework for control of linear finite-dimensional state-space systems. Specifically, there exists a controller $u=Kx$ such that $\dot x=Ax+Bu$ is stable if and only if there exists some $P>0$ and $Z$ such that $AP+PA^T+BZ+Z^TB^T<0$. This LMI follows directly from the dual version of the Lyapunov inequality $AP+PA^T<0$ via the variable substitution $Z=KP$ ($K$ is then given by $K=ZP^{-1}$). If $A(\delta)$ and $B(\delta)$ are uncertain, $\delta \in \Delta$, then we search for $P(\delta)$ and $Z(\alpha)$ (or a fixed $P$ for quadratic stability) and the inequality must hold for all $\delta \in \Delta$ - a problem which is more difficult, but still convex in the variables $P$ and $Z$. LMIs of this form were introduced in~\cite{bernussou_1989} and are the basis for a majority of LMI methods for controller synthesis (See Chapter 5 Notes in~\cite{boyd_book} for a discussion). The question, then, is how to obtain similar results for control of time-delay systems.

Our approach is to think of the delay system evolving on a Hilbert space $x \in X$ as
\[
\dot x(t)=\mathcal{A}x(t)+\mathcal{B}u(t).
\]
 We seek an operator $\mathcal{K}$ such that the feedback $u=\mathcal{K}x$ is stabilizing. Note that the input $u$ can also be infinite-dimensional so that we may represent systems with input delay in the same framework (using a Dirac operator for $\mathcal{B}$). We also note that we use full-state feedback, which assumes that measurements are retained for a period equal to the value of delay. Since such full-state measurements are often not available, ultimately the framework must include output feedback control - a more difficult problem.

In the Hilbert space framework, then, and focusing on the first terms, we seek a stabilizability condition of the form $\ip{x}{\mathcal{A}\mathcal{P}x}+\ip{\mathcal{A}\mathcal{P}x}{x}+\ip{x}{\mathcal{B}\mathcal{Z}x}+\ip{\mathcal{B}\mathcal{Z}x}{x}\le 0$ for all solutions $x\in X$. To create and test such an inequality, the first step is to establish and test a dual Lyapunov inequality of the form $\ip{x}{\mathcal{A}\mathcal{P}x}+\ip{\mathcal{A}\mathcal{P}x}{x}\le 0$ which guarantees stability of $\dot x=\mathcal{A}x$. Construction and testing of such a dual Lyapunov test is the main contribution of this paper. Due to space constraints, controller synthesis itself will be treated separately, but some early results on synthesis can be found in~\cite{peet_2009TDS}. In addition, while we discuss enforcement of the operator inequalities, this is not the main focus of the paper, which relies on restatement of existing results, primarily from~\cite{peet_2014ACC}. Indeed, we emphasize that the contribution of the paper (Theorems~\ref{thm:dual} and~\ref{thm:dual_MD}) is not a specific numerical method for determining stability of time-delay systems, but rather a new Lyapunov framework for solving the problem of controller synthesis. Moreover, the conditions are deliberately formulated in such a way that alternative approaches such as~\cite{Gu_1997},\cite{gouaisbaut_2009},\cite{seuret_2013} may also be applied in order to test stability and obtain stabilizing controllers.
Finally, we note that in abstract space, there have been a number of results on dual and adjoint systems~\cite{bensoussan_book}. Unfortunately, however, these dual systems are not delay-type systems and there is no clear relationship between stability of these adjoint and dual systems and stability of the original delayed system.
%

This paper is organized as follows. In Sections~\ref{sec:Lyapunov} and~\ref{sec:Framework} we develop a mathematical framework for expressing Lyapunov-based stability conditions as operator inequalities. In Section~\ref{sec:duality} we show that given additional constraints on the Lyapunov operator, satisfaction of the dual Lyapunov inequality $\ip{x}{\mathcal{A}\mathcal{P}x}+\ip{\mathcal{A}\mathcal{P}x}{x}\le0$ proves stability of the delayed system. In Sections~\ref{sec:structured_SD} and~\ref{sec:structured_MD} we define a restricted class of Lyapunov functionals and operators which are valid for the dual stability condition in both the single-delay and mutliple-delay cases. In Sections~\ref{sec:dual_stability_SD} and~\ref{sec:dual_stability_MD} we apply these classes of operators to obtain dual stability conditions. These conditions are formulated as positivity and negativity of Lyapunov functionals and may be considered the primary contribution of the paper. We also note that the dual stability conditions have a tridiagonal matrix structure which is distinct from standard Lyapunov-Krasovskii forms and may potentially be exploited to increase performance when studying systems with a large number of delays. In Sections~\ref{sec:positivity}, ~\ref{sec:positivity_LMI}, and~\ref{sec:spacing}, we show how SOS-based methods can be used to parameterize positive Lyapunov functionals and thereby enforce the inequality conditions in Sections~\ref{sec:dual_stability_SD} and~\ref{sec:dual_stability_MD}. Finally, in Sections~\ref{sec:dual_LMI_SD} and~\ref{sec:dual_LMI_MD}, we summarize our results with a set of LMI conditions for dual stability in both the single and multiple-delay cases. Section~\ref{sec:toolbox} describes our Matlab toolbox, available online, which facilitates construction and solution of the LMIs. Section~\ref{sec:validation} applies the results to a variety of stability problems and verifies that the dual stability test is not conservative.

\section{Notation}\label{sec:Notation}
Standard notation includes the Hilbert spaces $L_2^m[X]$ of square integrable functions from $X$ to $\R^m$  and $W^m_2[X]:=\{x\, :\, x, \dot x \in L_2^m[X]\}$. We use $L_2^m$ and $W_2^m$ when domains are clear from context. We also use the the extensions $L_2^{n \times m}[X]$ and $W_2^{n \times m}[X]$ for matrix-valued functions which map to $\R^{n\times m}$. $\mathcal{C}[X] \supset W_2[X]$ denotes the continuous functions on $X$. $S^n \subset \R^{n \times n}$ denotes the symmetric matrices. $I_n \in \S^n$ denotes the identity matrix. $0_{n\times m} \in \R^{n\times m}$ is the matrix of zeros with shorthand $0_{n}:=0_{n \times n}$. For a natural number, $N \in \N$, we adopt the index shorthand notation which denotes $[K]=\{1,\cdots,K\}$. Some additional notation is defined throughout the paper with a selected subset summarized in the Appendix.

%

\section{Lyapunov Stability of Time-Delay Systems}\label{sec:Lyapunov}

In this paper, we consider stability of linear discrete-delay systems of the form
\begin{align}
\dot{x} (t) &= A_0x(t) + \sum_{i =1}^K A_i x(t-\tau_i)\;&\text{ for all } \; t\ge0 ,\notag \\
x(t) &= \phi(t)\;&\text{ for all } \; t \in [-\tau_K,0]  \label{eqn:delay_eqn}
\end{align}
where $A_i \in \R^{n\times n}$, $\phi \in \mathcal{C}[-\tau_K,0]$, $K\in \N$ and for convenience $\tau_1 < \tau_2 < \cdots < \tau_K$. We associate with any solution $x$ and any time $t \ge 0$, the `state' of System~\eqref{eqn:delay_eqn}, $x_t \in \mathcal{C}[-\tau_K,0]$, where $x_t(s) = x(t+s)$. Although we only consider discrete-delay systems, the results of this paper may easily be extended to systems with distributed delay. For linear discrete-delay systems of the form~\eqref{eqn:delay_eqn}, the system has a unique solution for any $\phi \in \mathcal{C}[-\tau_K,0]$ and global, local, asymptotic and exponential stability are all equivalent.

Stability of Equations~\eqref{eqn:delay_eqn} may be certified through the use of Lyapunov-Krasovskii functionals - an extension of Lyapunov theory to systems with infinite-dimensional state-space. In particular, it is known that stability of linear time-delay systems is equivalent to the existence of a quadratic Lyapunov-Krasovskii functional of the form
\begin{align}
V(\phi) &= \int_{-\tau_K}^0\bmat{\phi(0)\\\phi(s)}^T M(s) \bmat{\phi(0)\\\phi(s)} ds \notag \\
 & +  \int_{-\tau_K}^0 \int_{-\tau_K}^0 \phi(s)^T N(s,\theta) \phi(\theta)\, ds\,d\theta, \label{eqn:complete_quad}
\end{align}
where the Lie (upper-Dini) derivative of the functional is negative along any solution $x$ of~\eqref{eqn:delay_eqn}. That is,
\[
\dot V(x_t) = \lim_{h \rightarrow 0}\frac{V(x_{t+h}) - V(x_t)}{h} \le 0
\]
for all $t \ge 0$.
Furthermore, the unknown functions $M$ and $N$ may be assumed to be continuous in their respective arguments everywhere except possibly at points $H := \{-\tau_1,\cdots,-\tau_K\}$.

For the dual stability conditions we propose in this paper, discontinuities in the unknown functions $M$ and $N$ pose challenges which make this form of Lyapunov-Krasovskii functional poorly suited to controller synthesis. For this reason, we use an alternative formulation of the necessary Lyapunov-Krasovskii functional better suited to the dual stability conditions we propose. Specifically, it has been shown~\cite{gu_2010} that existence of a positive decreasing Lyapunov-Krasovskii functional of the form in Eqn.~\eqref{eqn:complete_quad} implies the existence of a positive decreasing Lyapunov-Krasovskii functional of the form
\begin{align}
V(\phi) &= \tau_K \phi(0)^T  P \phi(0) +  \tau_K \sum_{i=1}^K \int_{-\tau_i}^0 \phi(0)^T Q_i(s)\phi(s) d s
+\tau_K\sum_{i=1}^K \int_{-\tau_i}^0 \phi(s)^T Q_i(s)^T \phi(0) ds \notag \\
 &+\tau_K\sum_{i=1}^K \int_{-\tau_i}^0 \phi_i(s)^T  S_i(s)\phi_i(s) + \sum_{i,j=1}^{K}\int_{-\tau_i}^0 \int_{-\tau_j}^0  \phi(s)^T R_{ij}(s,\theta)\phi(\theta)d \theta, \label{eqn:complete_quad2}
\end{align}
where the functions $Q_i$, $S_i$ and $R_{ij}$ may be assumed continuous on their respective domains of definition.

\section{A Mathematical Framework for Lyapunov Inequalities}\label{sec:Framework}
The use of Lyapunov-Krasovskii functionals can be simplified by considering stability in the semigroup framework - a generalization of the concept of differential equations. Although the results of this paper do not require the semigroup architecture, we adopt this notation in order to simplify the concepts and avoid unnecessary notation. Sometimes known as the `flow map', a `strongly continuous semigroup' is an operator, $S(t):Z \rightarrow Z$, defined by the Hilbert space $Z$, which represents the evolution of the state of the system so that for any solution $x$, $x_{t+s} = S(s) x_t$. Note that for a given $Z$, the semigroup may not exist even if the solution exists for any initial conditions in $Z$. Associated with a semigroup on $Z$ is an operator $\mathcal{A}$, called the `infinitesimal generator' which satisfies
\[
\frac{d}{dt}S(t)\phi = \mathcal{A} S(t)\phi
\]
for any $\phi \in X$. The space $X$ is often referred to as the domain of the generator $\mathcal{A}$, and is the space on which the generator is defined and need not be a closed subspace of $Z$. In this paper we will refer to $X$ as the `state-space'. For System~\eqref{eqn:delay_eqn}, following the approach in~\cite{curtain_book}, we define $Z_{m,n,K}:=\{\R^m \times L_2^{n}[-\tau_1,0]\times \cdots \times L_2^n[-\tau_K,0]\}$ and for $\{x,\phi_1,\cdots,\phi_K\}\in Z_{m,n,K}$, we define the following shorthand notation
\[
\bmat{x\\ \phi_i}:=\{x,\phi_1,\cdots,\phi_K\},
\]
which allows us to simplify expression of the inner product which we define to be
\[
\ip{\bmat{y\\ \psi_i}}{\bmat{x\\ \phi_i}}_{Z_{m,n,K}}=\tau_K y^T x + \sum_{i=1}^K \int_{-\tau_i}^0 \psi_i(s)^T\phi_i(s)ds.
\]
Furthermore, when $m=n$, we simplify the notation using $Z_{n,K}:=Z_{n,n,K}$. We may now conveniently write the state-space as
\[
X:=\left\{\bmat{x \\ \phi_i} \in Z_{n,K}\, : \, \phi_i \in W_2^n[-\tau_i,0] \text{ and } \phi_i(0)=x \text{ for all } i\in [K] \right\}.
\]
We furthermore extend this notation to say
\[
\bmat{x \\ \phi_i}(s)=\bmat{y \\ f(s,i)}
\]
if  $x=y$ and $\psi_i(s)=f(s,i)$ for $s \in [-\tau_i,0]$ and $i\in [K]$. This also allows us to compactly represent the infinitesimal generator, $\mathcal{A}$, of Eqn.~\eqref{eqn:delay_eqn} as
\[
\mathcal{A}\bmat{x\\ \phi_i}(s):= \bmat{A_0 x + \sum_{i=1}^K A_i \phi_i(-\tau_i)\vspace{2mm}\\ \dot \phi_i(s)}.
\]
Using these definitions of $\mathcal{A}$, $Z$ and $X$, for matrix $P$ and sufficiently smooth functions $Q_i,S_i,R_{ij}$, we define an operator $\mathcal{P}_{\{P,Q_i,S_i,R_{ij}\}}$ of the ``complete-quadratic'' type as
\begin{align}
&\left(\mathcal{P}_{\{P,Q_i,S_i,R_{ij}\}} \bmat{x\\ \phi_i}\right)(s) :=\bmat{  P x +  \sum_{i=1}^K \int_{-\tau_i}^0 Q_i(s)\phi_i(s) d s \\
\tau_K Q_i(s)^T x + \tau_K S_i(s)\phi_i(s) + \sum_{j=1}^K \int_{-\tau_j}^0 R_{ij}(s,\theta)\phi_j(\theta)\, d \theta }.
\end{align}
This notation allows us to associate $P,Q_i,S_i$ and $R_{ij}$ with the corresponding complete-quadratic functional in Eqn~\eqref{eqn:complete_quad2} as
\[
V(\phi)=\ip{\bmat{\phi(0)\\ \phi_i }}{\mathcal{P}_{\{P,Q_i,S_i,R_{ij}\}}\bmat{\phi(0)\\ \phi_i }}_{Z_{n,K}}.
\]

That is, the Lyapunov functional is defined by the operator $\mathcal{P}_{\{P,Q_i,S_i,R_{ij}\}}$ which is a variation of a classical combined multiplier and integral operator whose multipliers and kernel functions are given by $P,Q_i,S_i,R_{ij}$. The time-derivative of the complete-quadratic functional can similarly be represented using these operators as
\[
\dot V(\phi)=\ip{\bmat{\phi(0)\\ \phi_i }}{\mathcal{P}_{\{P,Q_i,S_i,R_{ij}\}}\mathcal{A} \bmat{\phi(0)\\ \phi_i }}_{Z_{n,K}}+\ip{\mathcal{A} \bmat{\phi(0)\\ \phi_i }}{\mathcal{P}_{\{P,Q_i,S_i,R_{ij}\}}  \bmat{\phi(0)\\ \phi_i }}_{Z_{n,K}}.
\]
The classical stability problem, then, states that the delay-differential Equation~\eqref{eqn:delay_eqn} is stable if there exists an $\alpha > 0$, matrix $P$ and functions $Q_i,S_i,R_{ij}$ such that $V(\phi)\ge \alpha \norm{\phi(0)}^2$ and $\dot V(\phi)\le 0$ for all $\phi\in \mathcal{C}[-\tau_K,0]^n$ such that $\bmat{\phi(0)\\ \phi_i} \in X$. In this paper, however, we seek to establish new stability conditions in a dual space - a problem which is formulated is formulated in the following Section.

%

\section{A Dual Stability Condition}\label{sec:duality}
Using the notation we have introduced in the proceeding section, we may compactly represent the dual stability condition which forms the main theoretical contribution of the paper.

\begin{thm} \label{thm:dual}
Suppose that $\mathcal{A}$ generates a strongly continuous semigroup on $Z$ with domain $X$. Further suppose there exists a bounded, positive and coercive linear operator $\mathcal{P} : X \rightarrow X$ which is self-adjoint with respect to the $Z$ inner product and satisfies
\[
\ip{\mathcal{A}\mathcal{P}z}{z}_Z+\ip{z}{\mathcal{A}\mathcal{P}z}_Z\le -\norm{z}^2_Z
\]
for all $z \in X$. Then a dynamical system which satisfies $\dot x(t) = \mathcal{A}x(t)$ generates an exponentially stable semigroup.
\end{thm}
\begin{IEEEproof}
Because $\mathcal{P}$ is coercive, bounded and self-adjoint, its inverse exists and is coercive, bounded and self-adjoint. Define the Lyapunov functional
\[
V(y) = \ip{y}{\mathcal{P}^{-1} y}\ge \alpha \norm{y}^2_Z
\]
which holds for some $\alpha>0$ and all $y\in X$. If $y(t)$ satisfies $\dot y(t)=\mathcal{A}y(t)$, then $V$ has time derivative
\begin{align}
\frac{d}{dt} V(y(t)) &= \ip{\dot y(t)}{\mathcal{P}^{-1}y(t)} + \ip{y(t)}{\mathcal{P}^{-1}\dot y(t)}\\
 &= \ip{\mathcal{A} y(t)}{\mathcal{P}^{-1}y(t)} + \ip{y(t)}{\mathcal{P}^{-1} \mathcal{A}y(t)}\\
 &= \ip{\mathcal{A} y(t)}{\mathcal{P}^{-1}y(t)} + \ip{\mathcal{P}^{-1}y(t)}{ \mathcal{A}y(t)}.
\end{align}
Now define $z(t)= \mathcal{P}^{-1} y(t) \in X$ for all $t\ge 0$. Then $y(t)=\mathcal{P}z(t)$ and since $\mathcal{P}$ is bounded and $\mathcal{P}^{-1}$ is coercive, there exist $\gamma,\delta>0$ such that
\begin{align}
\dot V(y(t)) &= \ip{\mathcal{A} y(t)}{\mathcal{P}^{-1}y(t)} + \ip{\mathcal{P}^{-1}y(t)}{ \mathcal{A}y(t)}\\
 &= \ip{\mathcal{A} \mathcal{P}z(t)}{z(t)} + \ip{z(t)}{ \mathcal{A}\mathcal{P}z(t)} \\
 &\le -\norm{z(t)}^2 \le -\frac{1}{\gamma}\ip{z(t)}{\mathcal{P}z(t)}=-\frac{1}{\gamma}\ip{y(t)}{\mathcal{P}^{-1}y(t)}\le -\frac{\delta}{\gamma} \norm{y(t)}^2.
\end{align}
Negativity of the derivative of the Lyapunov function implies exponential stability in the square norm of the state by, e.g.~\cite{curtain_book} or by the invariance principle.
\end{IEEEproof}

The advantage of the dual stability condition is that we replace $\mathcal{A}\mathcal{P}$ with $\mathcal{P}\mathcal{A}$. Although relatively subtle, this distinction allows convexification of the the controller synthesis problem. In the following section, we discuss how to parameterize operators which satisfy the conditions of Theorem~\ref{thm:dual}. We start with the constraints $\mathcal{P}=\mathcal{P^*}$ and $\mathcal{P}:X\rightarrow X$. Note that without significant restrictions on $P,Q_i,S_i,R_{ij}$, the operator $\mathcal{P}_{\{P,Q_i,S_i,R_{ij}\}}$ satisfies neither constraint.

\section{A Structured Operator: Single Delay}\label{sec:structured_SD}

In order to satisfy the conditions of the dual stability condition in Theorem~\ref{thm:dual}, we must restrict ourselves to a class of operators which are self-adjoint with respect to the inner-product defined on $Z_{n,K}$ and which preserve the structure of the state-space (map $X\rightarrow X$). We first consider the simpler case of a single delay. In this case, we have $Z_{m,n,1}=\R^m \times L_2^{n}$ with the $L_2^{m\times n}$ inner product and the state-space becomes $X:=\{\{x, \phi \} \in \R^n \times W_2^n[-\tau,0]\, : \,  \phi(0)=x   \}$ . To preserve the structure of $X$, we consider operators of the form

\begin{align}
&\left(\mathcal{P}\bmat{x \\ \phi}\right)(s):= \bmat{  \tau( R(0,0)+S(0))x +  \int_{-\tau}^0 R(0,s)\phi(s)d s \\ \tau R(s,0)\phi(0) + \tau S(s)\phi(s) + \int_{-\tau}^0 R(s,\theta)\phi(\theta)d \theta }\label{eqn:operator}
\end{align}
Clearly, we have that $\mathcal{P}$ is a bounded linear operator and if $S,R \in W_2^{n \times n}[-\tau,0]$ by inspection maps $\mathcal{P}: X \rightarrow X$. Furthermore, $\mathcal{P}$ is self-adjoint with respect to the $L_2^{2n}$ inner product, as indicated in the following lemma.

\begin{lem}\label{lem:selfadjoint}
Suppose $R(s,\theta) = R(\theta,s)^T$ and $S(s) \in \S^n$. Then the operator $\mathcal{P}$, as defined in Equation~\eqref{eqn:operator}, is self-adjoint with respect to the $L_2^{2n}$ inner product.
\end{lem}
\begin{IEEEproof}
The operator $\mathcal{P}: X \rightarrow X$ is self-adjoint with respect to the inner product $\ip{\cdot}{\cdot}_{L_2^{2n}}$ if
\[
\ip{\bmat{y \\\psi}}{\mathcal{P}\bmat{x \\ \phi}}_{L_2^{2n}}=\ip{\mathcal{P}\bmat{y \\\psi}}{\bmat{x \\ \phi}}_{L_2^{2n}}
\]
for any $\bmat{x \\ \phi},\bmat{y \\\psi} \in X$. By exploiting the structure of $\mathcal{P}$ and $X$, we have the following.{
\begin{align}
&\ip{\bmat{y \\\psi}}{\mathcal{P}\bmat{x \\ \phi}}_{L_2^{2n}}=\int_{-\tau}^0\bmat{y\\\psi(s)} \bmat{  \tau (R(0,0)+S(0))x +  \int_{-\tau}^0 R(0,\theta)\phi(\theta)d \theta \\
                \tau R(s,0)\phi(0) + \tau S(s)\phi(s) + \int_{-\tau}^0 R(s,\theta)\phi(\theta)d \theta }ds \\
&=\int_{-\tau}^0\bmat{y\\\psi(s)}
        \bmat{ \tau (R(0,0)+S(0)) & \tau R(0,s) \\ \tau R(s,0) & \tau S(s)}
                \bmat{x\\\phi(s)}ds + \int_{-\tau}^0\int_{-\tau}^0\bmat{y\\\psi(s)}
            \bmat{0&0\\0&R(s,\theta)}
                    \bmat{x\\\phi(\theta)}ds \,d\theta\\
&=\int_{-\tau}^0\left(\bmat{\tau (R(0,0) + S(0)) & \tau R(0,s) \\
                            \tau R(s,0) & \tau S(s)}^T
            \bmat{y\\\psi(s)}\right)^T  \bmat{x\\\phi(s)}ds + \int_{-\tau}^0\int_{-\tau}^0\bmat{0\\R(s,\theta)^T \psi(s)}^T \bmat{x\\\phi(\theta)}ds \, d\theta\\
&=\int_{-\tau}^0\bmat{\tau (R(0,0) + S(0) )y   + \tau R(0,s) \psi(s)\\
                         \tau R(s,0) y  + \tau S(s) \psi(s) + \int_{-\tau}^0R(\theta,s)^T \psi(\theta)d \theta }^T
            \bmat{x\\ \phi(s) } ds \\
&=\int_{-\tau}^0\bmat{\tau (R(0,0) + S(0) )y   + \int_{-\tau}^0 R(0,\theta) \psi(\theta)d\theta\\ \tau R(s,0) \psi(0) + \tau S(s) \psi(s)+\int_{-\tau}^0 R(s,\theta) \psi(\theta)d \theta}^T  \bmat{x\\\phi(s)}ds =\ip{\mathcal{P}\bmat{y \\\psi}}{\bmat{x \\ \phi}}_{L_2^{2n}}
\end{align}
}
\end{IEEEproof}
Note that the constraint that the operator be self adjoint significantly reduces the number of free variables. In the single delay case, we have made this explicit by replacing the variables $P$ and $Q$ with $P=\tau(R(0,0)+S(0))$ and $Q(s)=R(0,s)$. A natural question is whether the self-adjoint constraint introduces conservatism. While we cannot establish that the self-adjoint constraint is necessary and sufficient for stability, construction of converse Lyapunov functionals in, e.g.~\cite{gu_2010} indicate coupling between the functions and furthermore, the numerical results at the end of this paper indicate little if any conservatism in this constraint. We now apply this structured operator to Theorem~\ref{thm:dual} to obtain conditions on $S$ and $R$ for which stability holds.

\section{Dual Stability Conditions - Single Delay Case}\label{sec:dual_stability_SD}
In this section, we apply the structured operator in Section~\ref{sec:structured_SD} to the dual stability condition in Theorem~\ref{thm:dual} to establish conditions for stability in the single-delay case. Note that we do not yet discuss how to enforce these conditions. First recall that the generator, $\mathcal{A}$ is defined as
\begin{align}
\left(\mathcal{A}\bmat{x\\ \phi}\right)(s) &= \bmat{ A_0 x + A_1 \phi(t-\tau) \\ \frac{d}{d s} \phi(s)}.
\end{align}

\begin{thm}\label{thm:dual_SD}
Suppose there exist $\epsilon>0$ and functions $S\in W_2^{n\times n}[-\tau,0]$ and $R\in W_2^{n\times n}[[-\tau,0]\times [-\tau,0]]$ where $R(s,\theta) = R(\theta,s)^T$ and $S(s) \in \S^n$ such that $\ip{x}{\mathcal{P}x}_{L_2^{2n}} \ge\epsilon \norm{x}^2$ for all $x \in X$ and
\[
\ip{\bmat{x \\ \phi(-\tau) \\ \phi}}{\mathcal{D}\bmat{x \\ \phi(-\tau) \\ \phi}}_{L_2^{3n}}\le-\epsilon \norm{\bmat{x \\ \phi}}_{L_2^{2n}}^2
\] for all $\bmat{x\\ \phi} \in X$ where
\begin{align}
&\left(\mathcal{P}\bmat{x \\ \phi}\right)(s):= \bmat{  \tau( R(0,0)+S(0))x +  \int_{-\tau}^0 R(0,s)\phi(s)d s \\ \tau R(s,0)\phi(0) + \tau S(s)\phi(s) + \int_{-\tau}^0 R(s,\theta)\phi(\theta)d \theta }
\end{align}
and
\begin{align}
&\left(\mathcal{D}\bmat{x \\ y \\ \phi}\right)(s):= \bmat{ D_0 \bmat{x \\ y}  + \int_{-\tau}^0 V(s)\phi(s)ds\\
\tau V(s)^T\bmat{x & y}^T + \tau \dot S(s)\phi(s) + \int_{-\tau}^0 E(s,\theta)\phi(\theta)d \theta }
\end{align}
where
\begin{align}
&D_0:=\bmat{S_{11}+S_{11}^T & S_{12} \\
        S_{12}^T & S_{22}} ,\quad V(s)=\bmat{S_{13}(s)\\0},\\
&S_{11} := \tau A_0(R(0,0)+S(0)) + \tau A_1 R(-\tau,0) +\frac{1}{2} S(0), \notag \\
&S_{12} := \tau A_1 S(-\tau), \qquad S_{22} := - S(-\tau), \notag \\
&S_{13}(s) := A_0 R(0,s)+A_1 R(-\tau,s)+\dot R(s,0)^T, \notag \\
 &E(s,\theta):=\frac{d}{ds}R(s,\theta) + \frac{d}{d\theta} R(s,\theta). \notag
\end{align}
Then the system defined by Equation~\eqref{eqn:delay_eqn} is exponentially stable.
\end{thm}

\begin{IEEEproof}
Define the operators $\mathcal{A}$ and $\mathcal{P}$ as above. By assumption, the operator $\mathcal{P}$ is coercive. By Lemma~\ref{lem:selfadjoint}, $\mathcal{P}$ is self-adjoint and maps $X\rightarrow X$. This implies that by Theorem~\ref{thm:dual} the system is exponentially stable if
\[
\ip{\mathcal{A}\mathcal{P} \bmat{x\\ \phi}}{\bmat{x\\ \phi}}_{L_2^{2n}}+\ip{\bmat{x\\ \phi}}{\mathcal{A}\mathcal{P} \bmat{x\\ \phi}}_{L_2^{2n}} \le -\norm{\bmat{x\\ \phi}}_{L_2^{2n}}
\]
for all $\bmat{x\\ \phi} \in X$. We begin by constructing $\mathcal{A}\mathcal{P}x$.
\begin{align}
&\left(\mathcal{A}\mathcal{P}\bmat{x\\ \phi}\right)(s):= \bmat{y_1\\ y_2(s)} \notag \\
&y_1 =   \tau A_0( R(0,0)+S(0))x +  \int_{-\tau}^0 A_0 R(0,s)\phi(s)d s  \\
&\qquad \qquad + A_1\bbl(\tau R(-\tau,0)\phi(0) + \tau S(-\tau)\phi(-\tau)  + \int_{-\tau}^0 R(-\tau,\theta)\phi(\theta)d \theta\bbr)\\
&y_2(s)= \tau \frac{d}{ds}R(s,0)\phi(0) + \tau \dot S(s)\phi(s) +  \tau S(s)\dot \phi(s)+ \int_{-\tau}^0 \frac{d}{ds}R(s,\theta)\phi(\theta)d \theta.
\end{align}
Thus
\begin{align}
&\ip{\bmat{x\\ \phi}}{\mathcal{A}\mathcal{P}\bmat{x\\ \phi}} := \tau x^T y_1 + \int_{-\tau}^0\phi(s)^Ty_2(s)ds.
\end{align}
Examining these terms separately and using $x = \phi(0)$, we have {
\begin{align}
\tau x^T y_1 &=\tau^2 x^T A_0( R(0,0)+S(0))x   + \tau \int_{-\tau}^0 x^T A_0 R(0,s)\phi(s)d s  + \tau x^T A_1 \tau R(-\tau,0)\phi(0)\\
 & \qquad + \tau^2 x^T A_1 S(-\tau)\phi(-\tau) + \tau\int_{-\tau}^0 x^T A_1 R(-\tau,\theta)\phi(\theta)d \theta\\
&= \int_{-\tau}^0    \left( x^T \tau A_0( R(0,0)+S(0))x +  \tau  x^T A_0 R(0,s)\phi(s)\right)d s  \\
& + \int_{-\tau}^0 (x^T \tau A_1 R(-\tau,0)x + x^T \tau A_1 S(-\tau) \phi(-\tau) )ds + \int_{-\tau}^0  \tau x^T A_1 R(-\tau,s)\phi(s) ds\\
&=\int_{-\tau}^0 \bmat{x\\x(-\tau)\\x(s)}^T  \bmat{\tau A_0( R(0,0)+S(0)) + \tau A_1 R(-\tau,0) &  *^T &  *^T\\ \frac{\tau}{2} S(-\tau)^T A_1^T &0& *^T\\ \frac{\tau}{2}   R(0,s)^T A_0^T+\frac{\tau}{2}  R(-\tau,s)^T A_1^T&0&0 }
                    \bmat{x\\x(-\tau)\\x(s)}ds.
\end{align}
}
Examining the second term, we get {
\begin{align}
&\int_{-\tau}\phi(s)^Ty_2(s)ds=\int_{-\tau}\phi(s)^T\tau \left(  \dot R(s,0)\phi(0) + \dot S(s)\phi(s)\right)ds\\
&\hspace{2cm}+   \int_{-\tau}^0   \phi(s)^T \tau S(s)\dot \phi(s)ds+  \int_{-\tau}^0  \int_{-\tau}^0  \phi(s)^T \frac{d}{ds}R(s,\theta)\phi(\theta)d \theta ds\\
&=\int_{-\tau}\phi(s)^T\left(\tau \phi(s)^T \dot R(s,0)\phi(0) + \phi(s)^T \tau \dot S(s)\phi(s)\right)ds +  \frac{\tau}{2} x^T S(0) x - \frac{\tau}{2} \phi(-\tau)^T S(-\tau)\phi(-\tau)\\
 &\qquad- \half \int_{-\tau}^0 \phi(s)^T \tau \dot S(s) \phi(s)ds + \int_{-\tau}^0 \int_{-\tau}^0 \phi(s)^T \frac{d}{ds}R(s,\theta)\phi(\theta)d \theta ds\\
&= \int_{-\tau}^0  \bmat{x\\x(-\tau)\\x(s)}^T
            \hspace{-1.5mm}\bmat{\frac{1}{2} S(0) & *^T & *^T\\
             0& \hspace{-1.5mm}-\frac{1}{2} S(-\tau) &*^T\\
            \frac{\tau}{2} \dot R(s,0)&0& \hspace{-1.5mm}\frac{\tau}{2} \dot S(s) }
                     \bmat{x\\x(-\tau)\\x(s)}ds   + \int_{-\tau}^0 \int_{-\tau}^0 \phi(s)^T \frac{d}{ds}R(s,\theta)\phi(\theta)d \theta ds.
\end{align}
}
Combining both terms, and using symmetry of the inner product, we get
\begin{align}
&\ip{\bmat{x \\ \phi}}{\mathcal{A}\mathcal{P}\bmat{x \\ \phi}} +\ip{\mathcal{A}\mathcal{P}\bmat{x \\ \phi}}{\bmat{x \\ \phi}}=  \\
&\int_{-\tau}^0 \bmat{x\\x(-\tau)\\x(s)}^T \bmat{S_{11}+S_{11}^T  & S_{12} & \tau S_{13}(s)\\
        S_{12}^T & S_{22} & 0_n\\
        \tau S_{13}(s)^T & 0_n & \tau \dot S(s) }
\bmat{x\\x(-\tau)\\x(s)} ds
+ \int_{-\tau}^0 \int_{-\tau}^0 \phi(s)^T \left(\frac{d}{ds}R(s,\theta) + \frac{d}{d\theta} R(s,\theta)\right)\phi(\theta)d \theta ds\\
&=\ip{\bmat{x \\ \phi(-\tau) \\ \phi}}{\mathcal{D}\bmat{x \\ \phi(-\tau) \\ \phi}}_{L_2^{3n}}\le-\epsilon \norm{\bmat{x \\ \phi}}_{L_2^{2n}}^2.
\end{align}
Since $\ip{\mathcal{A}\mathcal{P} \bmat{x\\ \phi}}{\bmat{x\\ \phi}}_{L_2^{2n}}+\ip{\bmat{x\\ \phi}}{\mathcal{A}\mathcal{P} \bmat{x\\ \phi}}_{L_2^{2n}} \le -\norm{\bmat{x\\ \phi}}_{L_2^{2n}}$ for all $\bmat{x \\ \phi}\in X$ we conclude that the conditions of Theorem~\ref{thm:dual} are satisfied and hence System~\eqref{eqn:delay_eqn} is exponentially stable.
\end{IEEEproof}

\noindent \textbf{Dual Lyapunov-Krasovskii Form:} To summarize the results of Theorem~\ref{thm:dual_SD} in a more traditional Lyapunov-Krasovskii format, the system is stable if there exists a
\begin{align}
V(\phi) &= \int_{-\tau}^0 \bmat{\phi(0)\\ \phi(s)}^T \bmat{ \tau( R(0,0)+S(0))   & \tau R(0,s)\\
        \tau R(s,0) & \tau S(s) } \bmat{\phi(0)\\ \phi(s)} ds + \int_{-\tau}^0 \int_{-\tau}^0 \phi(s)^T R(s,\theta) \phi(\theta)d \theta ds
\end{align}
such that $V(\phi)\ge \norm{\bmat{\phi(0)\\ \phi}}^2$ and
{
\begin{align}
V_D(\phi)&=\int_{-\tau}^0 \bmat{\phi(0)\\ \phi(-\tau)\\ \phi(s)}^T \bmat{S_{11}+S_{11}^T  & S_{12} & \tau S_{13}(s)\\
        S_{12}^T & S_{22} & 0_n\\
        \tau S_{13}(s)^T & 0_n & \tau \dot S(s) }
\bmat{\phi(0)\\ \phi(-\tau)\\ \phi(s)} ds\\
&\hspace{2cm}+ \int_{-\tau}^0 \int_{-\tau}^0 \phi(s)^T \left(\frac{d}{ds}R(s,\theta) + \frac{d}{d\theta} R(s,\theta)\right)\phi(\theta)d \theta ds\le -\epsilon \norm{\bmat{\phi(0) \\ \phi}}.
\end{align}
}
Note that unlike the standard Lyapunov-Krasovskii functions, the derivative of the dual functional is tri-diagonal in both the single-delay and multiple-delay cases. When studying systems with a large number of delays, it may be possible to exploit this structure to offer performance improvement over the standard Lyapunov-Krasovskii form.

\section{A Structured Operator: Multiple Delay}\label{sec:structured_MD}
Now that we have considered the single delay case, we extend this result to multiple delays. In this case, the constraint that the operator be self-adjoint is expressed as a linear constraint on $P$ and the functions $Q_i$, $S_i$ and $R_{ij}$, none of which are eliminated as was done for the single delay case.
For the multiple delay case, recall the state-space is defined as
\[
X:=\left\{\bmat{x \\ \phi_i} \in Z_{n,K}\, : \, \phi_i \in W_2^n[-\tau_i,0] \text{ and } \phi_i(0)=x \text{ for all } i\in [K] \right\}.
\]
Likewise, recall the inner product on $Z$ for $x,y \in X$ as
\[
\ip{\bmat{y\\ \psi_i}}{\bmat{x\\ \phi_i}}_{Z_{m,n,K}}=\tau_K y^T x + \sum_{i=1}^K \int_{-\tau_i}^0 \psi_i(s)^T\phi_i(s)ds.
\]
Lastly, recall we consider operators of the form
\begin{align}
&\left(\mathcal{P}_{\{P,Q_i,S_i,R_{ij}\}} \bmat{x\\ \phi_i}\right)(s) :=\bmat{  P x +  \sum_{i=1}^K \int_{-\tau_i}^0 Q_i(s)\phi_i(s) d s \\
\tau_K Q_i(s)^T x + \tau_K S_i(s)\phi_i(s) + \sum_{j=1}^K \int_{-\tau_j}^0 R_{ij}(s,\theta)\phi_j(\theta)\, d \theta }.\label{eqn:operator_ndelay}
\end{align}


\begin{lem}\label{lem:selfadjoint_MD}
Suppose that $S_i\in W_2^{n\times n}[-\tau_i,0]$, $R_{ij}\in W_2^{n\times n}\left[[-\tau_i,0]\times[-\tau_j,0]\right]$ and $S_i(s)=S_i(s)^T$, $R_{ij}(s,\theta)=R_{ji}(\theta,s)^T$, $P=\tau_KQ_i(0)^T + \tau_KS_i(0)$ and $Q_j(s)=R_{ij}(0,s)$ for all $i,j\in [K]$. Then $\mathcal{P}_{\{P,Q_i,S_i,R_{ij}\}}$ is a bounded linear operator, maps $\mathcal{P}_{\{P,Q_i,S_i,R_{ij}\}}: X \rightarrow X$, and as defined in Equation~\eqref{eqn:operator_ndelay}, is self-adjoint with respect to the inner product defined on $Z_{n,K}$.
\end{lem}
\begin{IEEEproof}
To simplify the presentation, let $\mathcal{P}:=\mathcal{P}_{\{P,Q_i,S_i,R_{ij}\}}$. We first establish that $\mathcal{P}: X \rightarrow X$. If $\bmat{x\\\phi_i}\in X$, then $\phi_i \in \mathcal{C}[-\tau_i,0]$ and $\phi_i(0)=x$. Now if
\begin{align}
&\bmat{y\\ \psi_i(s)}=\left(\mathcal{P}\bmat{x\\ \phi_i}\right)(s)=\bmat{  P x +  \sum_{i=1}^K \int_{-\tau_i}^0 Q_i(s)\phi_i(s) d s \\
\tau_K Q_i(s)^T x + \tau_K S_i(s)\phi_i(s) + \sum_{j=1}^K \int_{-\tau_j}^0 R_{ij}(s,\theta)\phi_j(\theta)d \theta }
\end{align}
then since $P=\tau_K Q_i(0)^T + \tau_K S_i(0)$ and $Q_j(s)=R_{ij}(0,s)$, we have that
\begin{align}
\psi_i(0)&=\tau_K Q_i(0)^T x + \tau_KS_i(0)\phi_i(0) + \sum_{j=1}^K \int_{-\tau_j}^0 R_{ij}(0,\theta)\phi_j(\theta)d \theta\\
&=\left( \tau_K Q_i(0)^T + \tau_KS_i(0)\right)x + \sum_{j=1}^K \int_{-\tau_j}^0 R_{ij}(0,\theta)\phi_j(\theta)d \theta\\
&=P x + \sum_{j=1}^K \int_{-\tau_j}^0 Q_j(s)\phi_j(s)d s\\
&=y.
\end{align}
Since $S_i\in W_2^{n\times n}[-\tau_i,0]$, $R_{ij}\in W_2^{n\times n}\left[[-\tau_i,0]\times[-\tau_j,0]\right]$, $\phi_i\in W_2^n[-\tau_i,0]$, and hence we have $\bmat{y \\ \psi_i}\in X$ and hence $\mathcal{P}:X\rightarrow X$. Furthermore, boundedness of $Q_i$, $S_i$ and $R_{ij}$ implies boundedness of the linear operator $\mathcal{P}$.

Now, to prove that the operator $\mathcal{P}$ is self-adjoint with respect to the inner product $\ip{\cdot}{\cdot}_{Z_{n,K}}$, we show
\[
\ip{y}{\mathcal{P}x}_{Z_{n,K}}=\ip{\mathcal{P}y}{x}_{Z_{n,K}}
\]
for any $x,y \in X$. Using the properties $S_i(s)=S_i(s)^T$ and $R_{ij}(s,\theta)=R_{ji}(\theta,s)^T$, we have the following.

{
\begin{align}
&\ip{\bmat{y \\ \psi_i}}{\mathcal{P}\bmat{x \\ \phi_i}}_{Z_{n,K}}=\tau_K y^T\left(   P x +  \sum_{i=1}^K \int_{-\tau_i}^0 Q_i(\theta)\phi_i(\theta) d \theta\right)\\
   &\qquad + \sum_{i=1}^K \int_{-\tau_i}^0 \psi_i(s)\left(\tau_K Q_i(s)^T x + \tau_KS_i(s)\phi_i(s) + \sum_{j=1}^K \int_{-\tau_j}^0 R_{ij}(s,\theta)\phi_j(\theta)d \theta \right)\\
&=y^T \tau_K P x +  \sum_{i=1}^K \int_{-\tau_i}^0 y^T \tau_K Q_i(s)\phi_i(s) d s\\
   &+ \sum_{i=1}^K \int_{-\tau_i}^0 \left(\psi_i(s)^T \tau_K Q_i(s)^T x + \tau_K\psi_i(s) S_i(s)\phi_i(s) + \sum_{j=1}^K \int_{-\tau_j}^0 \psi_i(s)R_{ij}(s,\theta)\phi_j(\theta)d \theta \right)ds\\
&=\left(\tau_K Py + \sum_{i=1}^K \int_{-\tau_i}^0 \tau_K Q_i(s) \psi_i(s)ds  \right)^T  x \\
   &\qquad + \sum_{i=1}^K \int_{-\tau_i}^0 \left(y^T \tau_K Q_i(s)+ \tau_K\psi_i(s)^T S_i(s) + \sum_{j=1}^K \int_{-\tau_j}^0 \psi_j(\theta)^TR_{ji}(\theta,s)d\theta   \right)\phi_i(s) ds\\
&=\tau_K  \left(Py + \sum_{j=1}^K \int\limits_{-\tau_j}^0 Q_i(s) \psi_j(s)ds  \right)^T  x\\
 &\qquad \qquad +  \sum_{i=1}^K \int\limits_{-\tau_i}^0 \left( \tau_K Q_i(s)^T y + \tau_KS_i(s)^T \psi_i(s) +\sum_{j=1}^K \int\limits_{-\tau_j}^0 R_{ji}(\theta,s)^T\psi_j(\theta)d\theta \right)^T \phi_i(s) \,ds\\
&\qquad =\ip{\mathcal{P}\bmat{y\\\psi_i}}{\bmat{x\\ \phi_i}}_{Z_{n,K}}
\end{align}
}
\end{IEEEproof}

\section{The Dual Stability Condition for Multiple Delays}\label{sec:dual_stability_MD}
For the multiple-delay case, we apply the operator defined in Section~\ref{sec:structured_MD} to the dual stability condition in Theorem~\ref{thm:dual}. Here the generator, $\mathcal{A}$ is defined as
\begin{align}
\left(\mathcal{A} \bmat{x\\ \phi_i}\right)(s) &= \bmat{ A_0 x + \sum_{i=1}^k A_i \phi_i(-\tau_i) \\ \frac{d}{d s} \phi_i(s)}.
\end{align}

\begin{thm}\label{thm:dual_MD}
Suppose that there exist $S_i\in W_2^{n\times n}[-\tau_i,0]$ and $R_{ij}\in W_2^{n\times n}\left[[-\tau_i,0]\times[-\tau_j,0]\right]$ such that $S_i(s)=S_i(s)^T$ and $R_{ij}(s,\theta)=R_{ji}(\theta,s)^T$. Let $P=\tau_KQ_i(0)^T + \tau_KS_i(0)$ and $Q_j(s)=R_{ij}(0,s)$ for all $i,j\in [K]$. If $\ip{x}{\mathcal{P}x}_{Z_{n,K}} \ge \epsilon \norm{x}^2$ for all $x \in X$ and
\[
\ip{\bmat{\bmat{x \\ \phi_1(-\tau_1) \\ \vdots \\ \phi_k(-\tau_K)}\\ \phi_i}}{\mathcal{D}\bmat{\bmat{x \\ \phi_1(-\tau_1) \\ \vdots \\ \phi_k(-\tau_K)}\\ \phi_i}}_{Z_{{nK},n,K}}\le- \norm{\bmat{x \\ \phi_i}}_{Z_{n,K}}^2
\] for all $\bmat{x\\ \phi_i} \in X$ where
\begin{align}
\left(\mathcal{P} \bmat{x\\ \phi_i}\right)(s) =\bmat{  P x +  \sum_{i=1}^K \int_{-\tau_i}^0 Q_i(s)\phi_i(s) d s \\
\tau_K Q_i(s)^T x + \tau_K S_i(s)\phi_i(s) + \sum_{j=1}^K \int_{-\tau_j}^0 R_{ij}(s,\theta)\phi_j(\theta)d \theta }
\end{align}
and
\begin{align}
\mathcal{D}\bmat{\bmat{x \\ \phi_1(-\tau_1) \\ \vdots \\ \phi_k(-\tau_K)}\\ \phi_i}(s)=\bmat{\bmat{C_0 & C_{1} & \cdots &C_{k} \\C_{1}^T &-S_1(-\tau_1) & 0&0 \\
\vdots&0&\ddots&0 \\
C_{k^T} &0&0&-S_k(-\tau_K)}\bmat{x \\ \phi_1(-\tau_1) \\ \vdots \\ \phi_k(-\tau_K)} + \sum_{i=1}^K\int_{-\tau_i}^0 \bmat{B_{i}(s)\\ 0 \\ \vdots \\ 0}\phi_i(s)ds\vspace{3mm} \\
\tau_K B_{i}(s)^T x + \tau_K \dot S_i(s)\phi_i(s)+\sum_{j=1}^K \int_{-\tau_j}^0G_{ij}(s,\theta)\phi_j(\theta) d \theta }
\end{align}
%
where
\begin{align}
&C_0:= A_0 P + PA_0^T +\tau_K \sum_{i=1}^K (  A_i Q_i(-\tau_i)^T+ Q_i(-\tau_i)A_i^T + S_i(0)),   \notag\\
&C_{i}:=\tau_K A_iS_i(-\tau_i),\notag \\
&B_{i}(s):=A_0 Q_i(s) +\dot Q_i(s)+\sum_{j=1}^K R_{ji}(-\tau_j,s), \notag \\
&G_{ij}(s,\theta):=\frac{\partial}{\partial s}R_{ij}(s,\theta)+\frac{\partial}{\partial \theta}R_{ji}(s,\theta)^T,
\end{align}
then the system defined by Equation~\eqref{eqn:delay_eqn} is exponentially stable.
\end{thm}

\begin{IEEEproof}
Define the operators $\mathcal{A}$ and $\mathcal{P}$ as above. By Lemma~\ref{lem:selfadjoint_MD}, $\mathcal{P}$ is self-adjoint and maps $X \rightarrow X$. Since $\mathcal{P}$ is positive and coercive by assumption, this implies by Theorem~\ref{thm:dual} the system is exponentially stable if
\[
\ip{\mathcal{A}P \bmat{x\\ \phi_i}}{\bmat{x\\ \phi_i}}+\ip{\bmat{x\\ \phi_i}\\ \phi_i}{\mathcal{A}P \bmat{x\\ \phi_i}} \le -\norm{\bmat{x\\ \phi_i}}^2
\]
for all $\bmat{x\\ \phi_i} \in X$. We begin by constructing $(\mathcal{A}\mathcal{P}x)(s):= \bmat{y\\\psi_i(s)}$.
\begin{align}
&y =   A_0 P x +  \sum_{i=1}^K \int_{-\tau_i}^0 A_0 Q_i(s)\phi_i(s)d s  \\
&\qquad + \sum_{i=1}^K A_i\left(\tau_K Q_i(-\tau_i)^T x + \tau_K S_i(-\tau_i)\phi_i(-\tau_i) + \sum_{j=1}^K \int_{-\tau_j}^0 R_{ij}(-\tau_i,\theta)\phi_j(\theta)d \theta\right),\\
&\psi_i(s)= \tau_K \dot Q_i(s)^T x + \tau_K \dot S_i(s)\phi_i(s) +  \tau_K S_i(s)\dot \phi_i(s)+ \sum_{j=1}^K \int_{-\tau_j}^0 \frac{d}{ds}R_{ij}(s,\theta)\phi_j(\theta)d \theta.
\end{align}
Thus
\begin{align}
&\ip{ \bmat{x\\ \phi_i}}{\mathcal{A}\mathcal{P} \bmat{x\\ \phi_i}}:= \tau_K x^T y  + \sum_{i=1}^K \int_{-\tau_i}^0\phi_i(s)^T\psi_i(s)ds.
\end{align}
Examining these terms separately and using $x = \phi_i(0)$, we have {
\begin{align}
&x^T y =x^T A_0Px +\sum_{i=1}^K \int_{-\tau_i}^0 x^T A_0Q_i(s)\phi_i(s)ds  + \sum_{i=1}^K \tau_K  x^T A_i Q_i(-\tau_i)^Tx\\
  & \qquad + \sum_{i=1}^K \tau_K x^T A_i S_i(-\tau_i)\phi_i(-\tau_i)+ \sum_{i=1}^K \int_{-\tau_i}^0 \sum_{j=1}^K x^T A_j R_{ji}(-\tau_j,\theta)\phi_i(\theta)d \theta
\end{align}
}
Examining the second term, we get {
\begin{align}
&\sum_{i=1}^K \int_{-\tau_i}^0 \phi_i(s)^T \psi_i(s)ds\\
&=\sum_{i=1}^K \tau_K \int_{-\tau_i}^0 \phi_i(s)^T \dot Q_i(s)^T x \,ds
+\sum_{i=1}^K \tau_K \int_{-\tau_i}^0 \phi_i(s)^T \dot S_i(s)\phi_i(s)ds+\sum_{i=1}^K \tau_K \int_{-\tau_i}^0 \phi_i(s)^T S_i(s) \dot \phi_i(s)ds\\
& \qquad  +\sum_{i,j} \int_{-\tau_i}^0 \int_{-\tau_j}^0 \phi_i(s)^T \frac{\partial}{\partial s}R_{ij}(s,\theta)\phi_i(\theta)\,ds\,d\theta\\
&=\sum_{i=1}^K \tau_K \int_{-\tau_i}^0 \phi_i(s)^T \dot Q_i(s)^T x \,ds
+\frac{\tau_K}{2} \sum_{i=1}^K \int_{-\tau_i}^0 \phi_i(s)^T \dot S_i(s)\phi_i(s)ds+ \frac{\tau_K}{2} x^T \sum_{i=1}^K S_i(0) x\\
& \qquad -\frac{\tau_K}{2} \sum_{i=1}^K \phi_i(-\tau_i)^T S_i(-\tau_i) \phi_i(-\tau_i)+\sum_{i,j} \int_{-\tau_i}^0 \int_{-\tau_j}^0 \phi_i(s)^T \frac{\partial}{\partial s}R_{ij}(s,\theta)\phi_i(\theta)\,ds\,d\theta
\end{align}
}
Combining both terms,
{
\begin{align}
&\ip{\bmat{x\\ \phi_i}}{\mathcal{A}\mathcal{P}\bmat{x\\ \phi_i}}_{Z_{n,K}}=\tau_K x^T y +\sum_{i=1}^K \int_{-\tau_i}^0 \phi_i(s)^T \psi_i(s)ds\\
&=x^T \left(\tau_K A_0P + \sum_{i=1}^K \tau_K^2 A_i Q_i(-\tau_i)^T +\frac{\tau_K}{2} \sum_{i=1}^K S_i(0) \right)x \\
&+  \tau_K^2\sum_{i=1}^K x^T A_i S_i(-\tau_i)\phi_i(-\tau_i) -\frac{\tau_K}{2}\sum_{i=1}^K \phi_i(-\tau_i)^T S_i(-\tau_i) \phi_i(-\tau_i)\\
&+\tau_K\sum_{i=1}^K \int_{-\tau_i}^0 x^T \left(A_0Q_i(s) + \dot Q_i(s)+\sum_{j=1}^K A_j R_{ji}(-\tau_j,s) \right)\phi_i(s)ds \\
 &+\frac{\tau_K}{2} \sum_{i=1}^K \int_{-\tau_i}^0 \phi_i(s)^T \dot S_i(s)\phi_i(s)ds+\sum_{i,j} \int_{-\tau_i}^0 \int_{-\tau_j}^0 \phi_i(s)^T \frac{\partial}{\partial s}R_{ij}(s,\theta)\phi_i(\theta)\,ds\,d\theta
\end{align}
}

Combining this term with its adjoint, we recover
\begin{align}
&\ip{\mathcal{A}\mathcal{P}\bmat{x\\ \phi_i}}{\bmat{x\\ \phi_i}}_{Z_{n,K}}+\ip{\bmat{x\\ \phi_i}}{\mathcal{A}\mathcal{P}\bmat{x\\ \phi_i}}_{Z_{n,K}}\hspace{-.5cm}=\ip{\bmat{\bmat{x \\ \phi_1(-\tau_1) \\ \vdots \\ \phi_k(-\tau_K)}\\ \phi_i}}{\mathcal{D}\bmat{\bmat{x \\ \phi_1(-\tau_1) \\ \vdots \\ \phi_k(-\tau_K)}\\ \phi_i}}_{Z_{{nK},n,K}}\hspace{-1cm}\le- \norm{\bmat{x \\ \phi_i}}_{Z_{n,K}}^2.
\end{align}
We conclude that all conditions of Theorem~\ref{thm:dual} are satisfied and hence System~\eqref{eqn:delay_eqn} is stable.
\end{IEEEproof}
In the following sections, we will show how positivity of $\mathcal{P}$ and negativity of $\mathcal{D}$ can be enforced using SDP when the functions $S_i$ and $R_{ij}$ are polynomial.\\

\noindent \textbf{Dual Lyapunov-Krasovskii Form:} To summarize the results of Theorem~\ref{thm:dual_MD} in a more traditional Lyapunov-Krasovskii format, the system is stable if there exists a
\begin{align}
V(\phi) &= \tau_K \phi(0)^T  P \phi(0) +  \tau_K \sum_{i=1}^K \int_{-\tau_i}^0 \phi(0)^T Q_i(s)\phi(s) d s
+\tau_K\sum_{i=1}^K \int_{-\tau_i}^0 \phi(s)^T Q_i(s)^T \phi(0) ds \notag \\
 &+\tau_K\sum_{i=1}^K \int_{-\tau_i}^0 \phi_i(s)^T  S_i(s)\phi_i(s) + \sum_{i,j=1}^{K}\int_{-\tau_i}^0 \int_{-\tau_j}^0  \phi(s)^T R_{ij}(s,\theta)\phi(\theta)d \theta,
\end{align}
such that $V(\phi)\ge \norm{\bmat{\phi(0)\\ \phi_i}}^2$ and
{
\begin{align}
V_D(\phi)&=\tau_K \phi(0)^T C_0 \phi(0) +  2\tau_K\sum_{i=1}^K \phi(0)^T C_i \phi_i(-\tau_i) - \tau_K \sum_{i=1}^K \phi_i(-\tau_i)^T S_i(-\tau_i) \phi_i(-\tau_i)\\
&+2\tau_K\sum_{i=1}^K \int_{-\tau_i}^0 \phi(0)^T B_i(s) \phi_i(s)ds +\tau_K \sum_{i=1}^K \int_{-\tau_i}^0 \phi_i(s)^T \dot S_i(s)\phi_i(s)ds\\
&+\sum_{i,j} \int_{-\tau_i}^0 \int_{-\tau_j}^0 \phi_i(s)^T G_{ij}(s,\theta)\phi_i(\theta)\,ds\,d\theta\le -\norm{\bmat{\phi(0)\\ \phi_i}}^2.
\end{align}
}

\section{SOS Conditions for Positivity on $Z_{m,n,K}$}\label{sec:positivity}
In the proceeding two sections, we have shown that stability of the multiple delay system is implied by the existence of an operator $\mathcal{P}_{\{P,Q_i,S_i,R_{ij}\}}$, which is positive on $Z_{n,n,K}$ and such that $\mathcal{D}$ is negative definite on $Z_{nK,n,K}$
and where $\mathcal{D}$ has a structure similar to $\mathcal{P}$ and is defined by functions which are linear transformations of the functions $P,Q_i,S_i,R_{ij}$. The challenge, then, is to search for the functions $P,Q_i,S_i,R_{ij}$ such that $\mathcal{P}$ is positive and $\mathcal{D}$ is negative. In this section, we discuss how to enforce positivity of $\mathcal{P}$ by assuming $Q_i,S_i,R_{ij}$ are polynomials and defining constraints on the coefficients of these polynomials in a form expressible as a semidefinite program.

Roughly speaking, our approach is to use positive matrices to parameterize a cone of operators with a square root defined on the appropriate inner product. For example, in $L_2$, if $Q>0$ is a positive matrix, it has a square root and hence if we define $V(z)=\ip{z}{Qz}$, we have $V(z)=\ip{z}{Qz}=\ip{z}{P^TPz}=\ip{Pz}{Pz}\ge 0$. Hence $(\mathcal{P}z)(s):=Qz(s)$ defines a positive operator. For $L_2[X]$, we generalize this approach using more complicated vectors of operators to obtain forms such as $V(z)=\ip{\mathcal{Z}(z)}{Q\mathcal{Z}(z)}_{L_2}$, as will be discussed in the following sections.
Unfortunately, however, positivity in the inner product on $Z_{m,n,K}$ is difficult to enforce directly. The reason, through some abuse of notation, is that unlike the $L_2$ inner product, for an arbitrary matrix $P$, $\ip{z}{P^T P z}_{Z_{m,n,K}}\neq \ip{P z}{Pz}_{Z_{m,n,K}}$. This difficulty may be overcome, however, by defining a transformation from $Z_{m,n,K}$ to $\R^m \times L_2^n[-\tau_K,0]$. Hence, our positive operators on elements of $Z_{m,n,K}$ will be a combination of a transformation from $Z_{m,n,K}$ to $\R^m \times L_2^n[-\tau_K,0]$ and a positive quadratic form defined on the space $\R^m \times L_2^n[-\tau_K,0]$.

First, consider the operator, $\mathcal{P}:X\rightarrow X$,
\begin{align}
&\left(\mathcal{P} \bmat{x\\ \phi_i}\right)(s) =\bmat{  P x +  \sum_{i=1}^K \int_{-\tau_i}^0 Q_i(s)\phi_i(s) d s \\
\tau_K Q_i(s)^T x + \tau_K S_i(s)\phi_i(s) + \sum_{j=1}^K \int_{-\tau_j}^0 R_{ij}(s,\theta)\phi_j(\theta)d \theta }.
\end{align}

Then, for $\bmat{x\\ \phi_i} \in X$, if we have that $\phi_i(s)=\phi_j(s)=\phi(s)$ for all $i,j\in[K]$ and $s\in [-\tau_K,0]$, we have the obvious representation
\begin{align}
&\ip{\bmat{x\\ \phi_i(s)}}{\mathcal{P}\bmat{x\\ \phi_i(s)}}_{Z_{m,n,K}}\\
&=\int_{-\tau_K}^0\bmat{x \\ \phi(s)}^T M(s) \bmat{x \\ \phi(s)}ds + \int_{-\tau_K}^0\int_{-\tau_K}^0\phi(s)^T N(s,\theta)\phi(s)ds\,d\theta
\end{align}
where
\begin{align}
M(s)&=\begin{cases}
\bmat{P & \tau_K \sum_{j=i}^k Q_j(s)\\
\tau_K \sum_{j=i}^k Q_j(s)^T & \tau_K \sum_{j=i}^k S_j(s) } & s \in [-\tau_{i}, -\tau_{i-1}]\\
\end{cases}\\
N(s,\theta)&=\begin{cases}
\sum_{l=i}^k \sum_{m=j}^k R_{lm}(s,\theta) & s \in [-\tau_{i}, -\tau_{i-1}],\, \theta \in [-\tau_{j}, -\tau_{j-1}]\\
\end{cases}\\
\end{align}
Then if we constrain $M$ and $N$ to define a positive operator on $\R^m \times L_2^n[-\tau_K,0]$, $\mathcal{P}$ will define a positive operator on $Z_{m,n,K}$.

Unfortunately, while $\phi_i(s)=\phi_j(s)=\phi(s)$ holds for solutions of Eqn~\eqref{eqn:delay_eqn}, elements of the dual state $z=\mathcal{P}\phi$ does not necessarily satisfy this property. Indeed, for an arbitrary $\bmat{x & \phi_i}^T \in X$, the restriction $\phi_i(s)=\phi_j(s)$ would place unreasonable additional constraints on the variables $Q_i$, $S_i$ and $R_{ij}$. For this reason, we instead perform a change of variables to obtain
\begin{align}
&\ip{\bmat{x\\ \phi_i(s)}}{\mathcal{P}\bmat{x\\ \phi_i(s)}}_{Z_{m,n,K}}\\
&=\int_{-\tau_K}^0\bmat{x \\ \hat \phi(s)}^T M(s) \bmat{x \\ \hat \phi(s)}ds + \int_{-\tau_K}^0\int_{-\tau_K}^0\hat\phi(s)^T N(s,\theta)\hat \phi(s)ds\,d\theta
\end{align}
where if define $a_i=\frac{\tau_i-\tau_{i-1}}{\tau_i}$, then
\begin{align}
M(s)&=\begin{cases}
\bmat{P & \frac{\tau_K}{a_i} Q_i(\frac{s+\tau_{i-1}}{a_i})\\
\frac{\tau_K}{a_i}  Q_i(\frac{s+\tau_{i-1}}{a_i})^T & \frac{\tau_K}{a_i} S_i(\frac{s+\tau_{i-1}}{a_i}) } & s \in [-\tau_i, -\tau_{i-1}]\\
\end{cases}\\
N(s,\theta)&=\begin{cases}
R_{ij}(\frac{s+\tau_{i-1}}{a_i},\frac{\theta+\tau_{j-1}}{a_j}) & s \in [-\tau_i, -\tau_{i-1}],\, \theta \in [-\tau_j, -\tau_{j-1}]\\
\end{cases}\\
\end{align}
and
\[
\hat \phi(s)=\begin{cases}
\phi_i(\frac{s+\tau_{i-1}}{a_i}) & s \in [-\tau_i, -\tau_{i-1}].\\
\end{cases}
\]
Thus, if $M$ and $N$ define a positive operator on $\R^m \times L_2^n[-\tau_K,0]$, then $\mathcal{P}$ defines a positive operator on $Z_{m,n,K}$.
Indeed, it can be shown that positivity of the operator $\mathcal{P}_{\{P,Q_i,S_i,R_{ij}\}}$ on $Z_{m,n,K}$ is equivalent~\cite{gu_2010} to positivity of the multiplier and  integral operator defined by the piecewise-continuous functions $M$ and $N$ on $L_2^n[-\tau_K,0]$ where we assume the $\hat \phi_i$ are all independent. To simplify notation, we will denote the transformation between $P,Q_i,S_i,R_{ij}$ and $M,N$ as
\[
\{M,N\}:=\mathcal{L}_1(P,Q_i,S_i,R_{ij})
\]
if $a_i=\frac{\tau_i-\tau_{i-1}}{\tau_i}$ and
\begin{align}
M(s)&=\begin{cases}
\bmat{P & \frac{\tau_K}{a_i} Q_i(\frac{s+\tau_{i-1}}{a_i})\\
\frac{\tau_K}{a_i}  Q_i(\frac{s+\tau_{i-1}}{a_i})^T & \frac{\tau_K}{a_i} S_i(\frac{s+\tau_{i-1}}{a_i}) } & s \in [-\tau_i, -\tau_{i-1}]\\
\end{cases}\notag \\
N(s,\theta)&=\begin{cases}
R_{ij}(\frac{s+\tau_{i-1}}{a_i},\frac{\theta+\tau_{j-1}}{a_j}) & s \in [-\tau_i, -\tau_{i-1}],\, \theta \in [-\tau_j, -\tau_{j-1}]\\
\end{cases}\label{eqn:linop1}
\end{align}

\begin{lem}
Let $\{M,N\}:=\mathcal{L}_1(P,Q_i,S_i,R_{ij})$ and
\begin{equation}
\left(\mathcal{P}_{M,N}x\right)(s):= M(s)x(s) + \int_{-\tau_K}^0 \bmat{0_{n} & 0_n\\ 0_n & N(s,\theta)}x(\theta)d \theta.\vspace{-1mm}
\end{equation}
If $\ip{x}{\mathcal{P}_{M,N}x}_{L_2^{m+n}}\ge \alpha \norm{x}_{L_2^{m+n}}^2$ for some $\alpha>0$ and all $x \in \R^m \times L_2^n[-\tau_K,0]$,  then $\ip{x}{\mathcal{P}_{\{P,Q_i,S_i,R_{ij}\}}x}_{Z_{m,n,K}}\ge \alpha \norm{x}^2_{Z_{m,n,K}}$ for all $x \in Z_{m,n,K}$.
\end{lem}
\begin{IEEEproof}
The proof follows directly from the observation that $\norm{\bmat{x\\ \hat \phi}}_{L_2^{m+n}}^2 = \norm{\bmat{x\\ \phi_i}}_{Z_{m,n,K}}^2$.
\end{IEEEproof}
Note that if $Q_i, S_i$ and $R_{ij}$ are polynomials with variable coefficients, then the constraint $\{M,N\}=\mathcal{L}_1(P,Q_i,S_i,R_{ij})$ defines a linear equality constraint between the coefficients of $Q_i, S_i$ and $R_{ij}$ and the coefficients of the polynomials which define $M$ and $N$. In the following section, we will discuss how to enforce positivity of operators on $\R^m \times L_2^n[-\tau_K,0]$ defined by piecewise-polynomial multipliers and kernels.

\section{LMI conditions for Positivity of Multiplier and Integral Operators}\label{sec:positivity_LMI}
In this Section, we define LMI-based conditions for positivity of operators of the form
\begin{equation}
\left(\mathcal{P}_{M,N}x\right)(s):= M(s)x(s) + \int_{-\tau_K}^0 N(s,\theta)x(\theta)d \theta.\label{eqn:operator_simple}\vspace{-1mm}
\end{equation}
where $x \in L_2^n[-\tau_K,0]$ and $M$ and $N$ are continuous except possibly on $s,\theta\in \{-\tau_1,\cdots-\tau_K\}$. In the following, for square-integrable functions $M,N$, we will retain the slightly overloaded notation $\mathcal{P}_{M,N}$ as defined in Equation~\eqref{eqn:operator_simple}. Note that we initially consider positivity of the operator on $L_2^{m+n}[-\tau_K,0]$ and not the subspace $\R^m \times L_2^n[-\tau_K,0]$.

Our approach to positivity is based on the observation that a positive operator will always have a square root. If we assume that this square root is also of the form of operator~\eqref{eqn:operator_simple} with functions $M$ and $N$ piecewise-polynomial of bounded degree, then the results of this section give necessary and sufficient conditions for the positivity of~\eqref{eqn:operator_simple}. Note that although this assumption is restrictive, it is unclear whether it implies conservatism. For example, while not all positive polynomials are Sum-of-Squares, any positive polynomial can be approximated arbitrarily well in the sup norm on a bounded domain by a polynomial with a polynomial ``root''.

\begin{thm}\label{thm:pos_op_joint}
For any functions $Y_1: [-\tau_K,0] \rightarrow \R^{m_1 \times n}$ and $Y_2: [-\tau_K,0] \times [-\tau_K,0] \rightarrow \R^{m_2 \times n}$, square integrable on $[-\tau_K,0]$ with $g(s)\ge 0$ for $s \in [-\tau_K,0]$, suppose that \vspace{-1mm}
\begin{align}
M(s) &= g(s) Y_1(s)^T Q_{11} Y_1(s) \label{def:M} \\
N(s,\theta) &= g(s) Y_1(s)Q_{12}Y_2(s,\theta) + g(\theta)Y_2(\theta,s)^T Q_{12}^T Y_1(\theta) \notag \\
&\qquad \qquad  + \int_{-\tau_K}^0 g(\omega)Y_2(\omega,s)^T Q_{22}Y_2(\omega,\theta) \, d\omega \label{def:N}\vspace{-1mm}
\end{align}
where $Q_{ij} \in \R^{m_i \times m_j}$ and
\[
Q=\bmat{Q_{11} & Q_{12}\\Q_{12}^T& Q_{22}} \ge 0.\vspace{-1mm}
\]
Then for $\mathcal{P}_{M,N}$ as defined in Equation~\eqref{eqn:operator_simple}, $\ip{x}{\mathcal{P}_{M,N}x}_{L_2^n} \ge 0$ for all $x \in L_2^n[-\tau_K,0]$.
\end{thm}
The proof of Theorem~\ref{thm:pos_op_joint} can be found in~\cite{peet_2014ACC}.

Theorem~\ref{thm:pos_op_joint} gives a linear parametrization of a cone of positive operators using positive semidefinite matrices. Note that there are few constraints on the functions $Y_1$ and $Y_2$. These functions serve as the basis for the multipliers and kernels found in the square root of $\mathcal{P}_{M,N}$. The class of multipliers and kernels defined by Theorem~\ref{thm:pos_op_joint} is thus determined by $Y_1$ and $Y_2$.

We now consider certain choices of $Y_1$ and $Y_2$ which yield piecewise-polynomials functions $M$ and $N$.\vspace{-1mm}

\subsection{Piecewise-Polynomials Multipliers and Kernels}\label{subsec:positivity_PC}

To define multipliers and kernels with discontinuities at known points, we divide the region of integration $[-\tau_K,0]$ into almost disjoint subregions $[-\tau_{i}, -\tau_{i-1}]$, $i \in [K]$ on which continuity holds and assume the functions are polynomial on these subregions. To do this, we introduce the indicator functions (not to be confused with the identity matrix)\vspace{-1mm}
\[
I_i(t) = \begin{cases}1 & t \in [-\tau_{i}, -\tau_{i-1}]\\ 0& \text{otherwise,} \end{cases} \quad i \in [K]\vspace{-2mm}
\]
and the vector of indicator functions $J = \bmat{I_1 & \cdots & I_K}^T$. We can now define the basis vectors $Y_1$ and $Y_2$ which define the positivity conditions in Theorem~\ref{thm:pos_op_joint}.
\[
Y_{1pc}(s)  =   Y_{1p}(s) \otimes J(s),\; Y_{2pc}(s,\theta)  =  Y_{2p}(s,\theta) \otimes J(s) \otimes J(\theta)\vspace{-1mm}
\]
where
\begin{equation}
Y_{1p}(s) = Y_d(s) \otimes I_n,\qquad Y_{2p}(s,\theta) = Y_d(s,\theta) \otimes I_n.\vspace{-1mm}
\end{equation}
and $Y_d(s)$ is a vector whose elements form a basis for the polynomials in variables $s$ of degree $d$ or less. e.g. The vector of monomials. Note for $s \in \R$, $Y_{d}:[-\tau_K,0]\rightarrow \R^{d+1}$, hence $Y_{1p}:[-\tau_K,0]\rightarrow \R^{n(d+1)\times n}$, and $Y_{1pc}:[-\tau_K,0]\rightarrow \R^{nK(d+1)\times n}$. Similarly, $Z_d(s,\theta)\in \R^{q}$ where $q=(d+1)(d+2)/2$, $Y_{2p}(s,\theta)\in \R^{nq\times n}$, and $Y_{2pc}(s,\theta)\in \R^{nKq\times n}$

\begin{thm}\label{thm:positivity_PC}
If $Y_{1}(s)=Y_{1pc}(s)$ and $Y_{2}(s,\theta) = Y_{2pc}(s,\theta)$ and $M$ and $N$ are defined as in Equations~\eqref{def:M} and~\eqref{def:N}, then $M$ and $N$ are piecewise-polynomial matrices ($\R^{n \times n}$) of degree $2d$ with possible discontinuities at $s,\theta \in \{-\tau_i\}_i$. In this case, if $g_i(s)\ge 0$ for $s \in [-\tau_i,-\tau_{i-1}]$, the functions $M$ and $N$ can be defined piecewise as\vspace{-1mm}
\[
M(s) = \begin{cases}M_i(s) & s \in [-\tau_{i}, -\tau_{i-1}]\end{cases}\vspace{-2mm}
\]
where\vspace{-1mm}
\[
M_{i} = g_i(s) Y_d(s)^T Q_{11,ii} Y_d(s)\vspace{-1.2mm}
\]
where $Q_{11,i,j}\in \R^{n(d+1)\times n(d+1)}$ is the $i,j$th block of $Q_{11}\in \S^{n(d+1)K}$. Likewise,\vspace{-1mm}
\[
N(s,\theta) = \begin{cases}N_{ij}(s,\theta) & s \in [-\tau_{i}, -\tau_{i-1}] \,\,\text{and} \,\,\theta \in [-\tau_{j}, -\tau_{j-1}]\end{cases}\vspace{-2mm}
\]
where\vspace{-1mm}
\begin{align}
&N_{ij} = g_i(s) Y_{1p}(s)Q_{12,i,(i-1)K+j}Y_{1p}(s,\theta) \\
&\qquad \qquad + g_j(\theta)Y_{2p}(\theta,s)^T Q_{12,(j-1)K+i,j}^T Y_{1p}(\theta)\\
& +  \sum_{l=1}^K \int_{-\tau_l}^{-\tau_{l-1}}  g_l(\omega) Y_{2p}(\omega_l,s)^T Q_{22,i+(l-1)K,j+(l-1)K}Y_{2p}(\omega_l,\theta) \, d\omega_l\vspace{-1mm}
\end{align}
where $Q_{12,i,j}\in \R^{n(d+1)\times nq}$ is the $i,j$th block of $Q_{12}\in \R^{n(d+1)K\times nqK}$ and $Q_{22,i,j}\in \R^{nq\times nq}$ is the $i,j$th block of $Q_{22}\in \S^{nqK}$.
\end{thm}
The proof of Theorem~\ref{thm:positivity_PC} can be found in~\cite{peet_2014ACC}.

For the intervals $s \in [-\tau_i,-\tau_{i-1}]$, the choice of $g_i$ is typically either $g_i(s)=1$ or $g_i=-(s+\tau_i)(s+\tau_{i-1})$. Inclusion of $g \neq 1$ is a variation of the classical Positivstellensatz approach to local positivity, as can be found in, e.g.~\cite{stengle_1973,schmudgen_1991,putinar_1993}. To improve accuracy, we typically use a combination of both although we may set $Q_{12},Q_{21},Q_{22}=0$ for the latter to reduce the number of variables. To simplify notation, throughout the remainder of the paper, we will use the notation $\{M,N\}\in \Xi_{d,n,K}$ to denote the LMI constraints on the coefficients of the polynomials $M,N$ implied by the conditions of Theorem~\ref{thm:positivity_PC} using both $g_i(s)=1$ and $g_i=-(s+\tau_i)(s+\tau_{i-1})$ as
\[
\Xi_{d,n,K}:=\{\{M,N\}\,:\,  \substack{ M=M_1+M_2,\, N=N_1+N_2,\, \text{ where $\{M_1,N_1\}$ and $\{M_2,N_2\}$ satisfy the}\\
\text{conditions of Thm.~\ref{thm:positivity_PC} with $g_i=1$ and $g_i=-(s+\tau_i)(s+\tau_{i-1})$, respectively.}}  \}
\]


\section{Spacing Functions and Mixed State-Space}\label{sec:spacing}\vspace{-1mm}

The result in Theorem~\ref{thm:pos_op_joint} as stated is a parametrization of operators which are positive on the space $L_2^n[-\tau_K,0]$. However, as in Section~\ref{sec:positivity}, we instead need to enforce positivity on the subspace $\R^m \times L_2^n[-\tau_K,0] \subset L_2^{m+n}[-\tau_K,0]$.

To enforce positivity on a subspace $X\subset L_2^n[-\tau_K,0]$, we turn to so-called ``spacing functions'' - a concept closely tied to projection operators.

\begin{thm}\label{prop:spacing}
Suppose $X$ is a closed subspace of a Hilbert space $Z$. Then $\ip{u}{\mathcal{P}u}\ge 0$ for all $u \in X$ if and only if there exist operators $\mathcal{M}$ and $\mathcal{T}$ such that $\mathcal{P}= \mathcal{M}+\mathcal{T}$ and $\ip{u}{\mathcal{M}u}\ge 0$ for all $u \in Z$ and $\ip{u}{\mathcal{T}u}=0$ for all $u \in X$.
\end{thm}
The proof of Theorem~\ref{thm:positivity_PC} can be found in~\cite{peet_2014ACC}.


This proposition implies that the class of operators which are positive on $X$ is the direct sum of the cone of operators, $\mathcal{M}$ which are positive on $Z$ and the space of operators, $\mathcal{T}$, which are orthogonal to $X$. Taking $Z=L_2^{n+m}$, we already know how to parameterize $\mathcal{M}$. The question, then, is how to parameterize the ``spacing'' operators $\mathcal{T}$.

\subsection{A Class of Spacing Functions}

For both the single-delay and multi-delay case, we enforce positivity on a subspace of the form $\R^m \times L_2^{n}[-\tau_K,0]\subset L_2^{m+n}[-\tau_K,0]$. For this subspace, we define a class of spacing functions as follows.
\begin{thm}\label{thm:mixed_pos}
Suppose that $F$ and $H$ are defined as
\begin{align}
&F(s) = \bmat{K(s) + \frac{1}{\tau_K}\int_{-\tau_K}^0 \int_{-\tau_K}^0 L_{11}(\omega,t)d\omega dt& \int_{-\tau_K}^0 L_{12}(\omega,s)d\omega \vspace{2mm} \\ \int_{-\tau_K}^0 L_{21}(s,\omega)d\omega & 0}\\
&H(s,\theta) =  - \bmat{L_{11}(s,\theta)&L_{12}(s,\theta)\\L_{21}(s,\theta)&0}
\end{align}
for some square-integrable functions $K$ and $L_{ij}$ where $K(s)\in \R^{m\times m}$,  $L_{11}(s,\theta)\in \R^{m\times m}$, and $L_{12}(s,\theta)\in \R^{m\times n}$ such that $\int_{-\tau_K}^0 K(s) ds =0$. Then if
\[
\mathcal{T}z(s):=F(s)z(s)+\int_{-\tau_K}^0H(s,\theta)z(\theta)\, d\theta
\]
then for any $z\in \R^m \times L_{2}^{n}$,
\[
\ip{z}{\mathcal{T}z}_{L_2^{m+n}}=0
\]
\end{thm}

\begin{IEEEproof}
The proof is straightforward. For $z(s)= \bmat{c & y(s)}^T$ with $c \in \R^m$ and $y \in L_2^{n}[-\tau_K,0]$, we have \vspace{-1mm}
\begin{align}
&\ip{z}{\mathcal{T}z}_{L_2^{m+n}}=\int_{-\tau_K}^0 \bmat{c \\ y(s)}^T \bmat{K(s) + \frac{1}{\tau_K}\int_{-\tau_K}^0 \int_{-\tau_K}^0 L_{11}(\omega,t)d\omega dt& \int_{-\tau_K}^0 L_{12}(\omega,s)d\omega \vspace{2mm} \\ \int_{-\tau_K}^0 L_{21}(s,\omega)d\omega & 0} \bmat{c \\ y(s)} ds \\
&\qquad \qquad \qquad  -\int_{-\tau_K}^0 \int_{-\tau_K}^0 \bmat{c \\ y(s)}^T \bmat{L_{11}(s,\theta)&L_{12}(s,\theta)\\L_{21}(s,\theta)&0} \bmat{c \\ y(\theta)} d\theta ds\\
&=\int_{-\tau_K}^0 \bmat{c \\ y(s)}^T \bmat{\frac{1}{\tau_K}\int_{-\tau_K}^0 K(\omega)d\omega & 0 \vspace{2mm} \\ 0 & 0} \bmat{c \\ y(s)} ds \\
&\qquad \qquad \qquad  +\int_{-\tau_K}^0 \int_{-\tau_K}^0 \bmat{c \\ y(s)}^T \bmat{L_{11}(s,\theta)-L_{11}(s,\theta)&L_{12}(s,\theta)-L_{12}(s,\theta)\\L_{21}(s,\theta)-L_{21}(s,\theta)&0} \bmat{c \\ y(\theta)} d\theta \,ds =0. \vspace{-1mm}
\end{align}
\end{IEEEproof}
For simplicity, we use $\{F,H\}\in \Theta_{m,n,K}$ to denote the conditions of Theorem~\ref{thm:mixed_pos} which, if $K$ and $L_{ij}$ are piecewise-polynomial matrices, is a set of linear equality constraints on the coefficients of the polynomials which define $F$ and $J$.
\[
\Theta_{m,n,K}:=\{\{F,H\}\,:\, F,H \text{ satisfy the conditions of Thm.~\ref{thm:mixed_pos}.}\}
\]
For convenience, given $\{F,H\} \in \Theta_{m,n,K}$, we define the operator $\mathcal{T}_{F,H}: L_2^{m,n}\rightarrow L_2^{m,n}$ as
\[
\left(\mathcal{T}_{F,H}z\right)(s):=F(s)z(s)+\int_{-\tau_K}^0H(s,\theta)z(\theta)\, d\theta.
\]

\section{SOS Conditions for Dual Stability in the Case of a Single Delay}\label{sec:dual_LMI_SD}
We now state an LMI representation of the dual stability condition for a single delay ($\tau=\tau_K$).
\begin{thm}\label{thm:dualLMI_SD}
Suppose there exist $d \in \N$, constant $\epsilon>0$, functions $S\in W_2^{n\times n}[-\tau,0]$, $R\in W_2^{n\times n}[[-\tau,0]\times [-\tau,0]]$, $\{F_1,H_1\}\in \Theta_{n,n,1}$, and $\{F_2,H_2\}\in \Theta_{2n,n,1}$ where $R(s,\theta) = R(\theta,s)^T$ and $S(s) \in \S^n$ such that
\[
\left\{M,N\right\} \in \Xi_{d,2n,1}
\]
and
\[
\left\{-D,-E\right\}\in \Xi_{d,3n,1}
\]
where
\begin{align}
&M(s)=\bmat{\tau R(0,0) + \tau S(0)& \tau R(0,s) \\ \tau R(s,0)&\tau S(s)}+F_1(s) -\epsilon I_{2n}, \\
&N(s,\theta)=\bmat{0_{n} & 0_{n}\\0_{n}&R(s,\theta)}+H_1(s,\theta),
\end{align}
\begin{align}
&D(s):=\bmat{D_0 & \tau V(s)\\
        \tau V(s)^T & \tau \dot S(s)+\epsilon I_n }+F_2(s) , \\
&D_0:=\bmat{S_{11}+S_{11}^T +\epsilon I_n& S_{12} \\
        S_{12}^T & S_{22}} ,\quad V(s)=\bmat{S_{13}(s)\\0},\\
&S_{11} := \tau A_0(R(0,0)+S(0)) + \tau A_1 R(-\tau,0) +\frac{1}{2} S(0),  \\
&S_{12} := \tau A_1 S(-\tau), \quad S_{22} := - S(-\tau),  \\
&S_{13}(s) := A_0 R(0,s)+  A_1 R(-\tau,s)+ \dot R(s,0)^T,  \\
 &E(s,\theta):=\bmat{0_{2n} & 0_{2n,n}\\0_{n,2n}&G(s,\theta)}+H_2(s,\theta) \\
 &G(s,\theta):=\frac{d}{ds}R(s,\theta) + \frac{d}{d\theta} R(s,\theta).
 \end{align}

Then the system defined by Equation~\eqref{eqn:delay_eqn} is exponentially stable.
\end{thm}
\begin{IEEEproof} Consider the operator
\begin{align}
&\left(\mathcal{P}\bmat{x \\ \phi}\right)(s):= \bmat{  \tau( R(0,0)+S(0))x +  \int_{-\tau}^0 R(0,s)\phi(s)d s \\ \tau R(s,0)\phi(0) + \tau S(s)\phi(s) + \int_{-\tau}^0 R(s,\theta)\phi(\theta)d \theta }
\end{align}
Since $\{F_1,H_1\}\in \Theta_{n,n,1}$, and  $\left\{M,N\right\} \in \Xi_{d,2n,1}$, by Lemma~\ref{prop:spacing} and Theorem~\ref{thm:positivity_PC}, we have for $x \in Z_{n,1}$
\begin{align}
\ip{x}{\mathcal{P}x}_{L_2^{2n}}-\epsilon \norm{x}^2 =\ip{x}{\left(\mathcal{P}+\mathcal{T}_{F_1,H_1}\right) x}_{L_2^{2n}}-\epsilon \norm{x}^2 =\ip{x}{\mathcal{P}_{M,N} x}_{L_2^{2n}}\ge 0.
\end{align}
This establishes that $\ip{x}{\mathcal{P}x}_{L_2^{2n}} \ge\epsilon \norm{x}^2$ for all $x \in X$. Similarly, examine the operator
\begin{align}
&\left(\mathcal{D}\bmat{x \\ y \\ \phi}\right)(s):= \bmat{ D_0 \bmat{x \\ y}  + \int_{-\tau}^0 V(s)\phi(s)ds\\
\tau V(s)^T\bmat{x & y}^T + \tau \dot S(s)\phi(s) + \int_{-\tau}^0 G(s,\theta)\phi(\theta)d \theta }.
\end{align}
Since $\{F_2,H_2\}\in \Theta_{2n,n,1}$, and  $\left\{-D,-E\right\} \in \Xi_{d,3n,1}$, we have for $x \in Z_{2n,n,1}$
\begin{align}
&\ip{\bmat{x \\ \phi(-\tau) \\ \phi}}{\mathcal{D}\bmat{x \\ \phi(-\tau) \\ \phi}}_{L_2^{3n}}+\epsilon \norm{\bmat{x \\ \phi}}^2 =\ip{\bmat{x \\ \phi(-\tau) \\ \phi}}{\left(\mathcal{D}+\mathcal{T}_{F_2,H_2}\right) \bmat{x \\ \phi(-\tau) \\ \phi}}_{L_2^{3n}}+\epsilon \norm{\bmat{x \\ \phi}}^2\\
& =\ip{\bmat{x \\ \phi(-\tau) \\ \phi}}{\mathcal{P}_{D,E} \bmat{x \\ \phi(-\tau) \\ \phi}}_{L_2^{3n}}\le 0.
\end{align}
This likewise establishes that
\[
\ip{\bmat{x \\ \phi(-\tau) \\ \phi}}{\mathcal{D}\bmat{x \\ \phi(-\tau) \\ \phi}}_{L_2^{3n}}\le-\epsilon \norm{\bmat{x \\ \phi}}^2
\]
for all $\bmat{x\\ \phi} \in X$. By assumption, $R(s,\theta) = R(\theta,s)^T$ and $S(s) \in \S^n$ and hence Theorem~\ref{thm:dual_SD} establishes exponential stability of Equation~\eqref{eqn:delay_eqn}.
\end{IEEEproof}

\section{SOS Conditions for Dual Stability in the Case of Multiple Delays}\label{sec:dual_LMI_MD}

\begin{thm}\label{thm:dualLMI_MD}
Suppose there exist $d \in \N$, constant $\epsilon>0$, matrix $P\in \R^{n\times n}$, functions $S_i, Q_i \in W_2^{n\times n}[-\tau_i,0]$, $R_{ij}\in W_2^{n\times n}\left[[-\tau_i,0]\times[-\tau_j,0]\right]$ for $i,j \in [K]$, $\{F_1,H_1\}\in \Theta_{n,n,K}$, and $\{F_2,H_2\}\in \Theta_{n(K+1),n,K}$ such that
\[
\left\{M,N\right\} \in \Xi_{d,2n,K} \qquad \text{and}\qquad \left\{-D,-E\right\}\in \Xi_{d,n(K+1),K},
\]
where
\begin{align}
M(s)&=M_0(s)+F_1(s)-\epsilon I_{2n} \quad \text{and} \quad  N(s,\theta)=\bmat{0_{n} & 0_{n}\\0_{n}&N_0(s,\theta)}+H_1(s,\theta),
\end{align}
where
\[
\{M_0,N_0\}:=\mathcal{L}_1(P-\epsilon I_n,Q_i,S_i-\epsilon I_n,R_{ij})
\]
and
\begin{align}
D(s)&=D_0(s)+F_2(s),\qquad E(s,\theta)=\bmat{0_{n(K+1)} & 0_{n(K+1),n}\\0_{n,n(K+1)}&E_0(s,\theta)}+H_2(s,\theta),
\end{align}
where
\[
\{D_0,E_0\}:=\mathcal{L}_1(D_1,V_i,\dot S_i + \epsilon I_n ,G_{ij})
\]

and where
\begin{align}
&D_1:=\bmat{C_{0}+C_0^T+\epsilon I_n & C_{1} & \cdots &C_{k} \\C_{1}^T &-S_1(-\tau_1) & 0&0 \\
\vdots&0&\ddots&0 \\
C_{k}^T &0&0&-S_k(-\tau_K)},\\
&C_{0}:= A_0 P +\tau_K \sum_{i=1}^K (  A_i Q_i(-\tau_i)^T + \half S_i(0)),   \notag\\
&C_{i}:=\tau_K A_iS_i(-\tau_i)\qquad i\in [K] \notag \\
&V_i(s):=\bmat{B_i(s)^T & 0 &\cdots &0}^T \qquad i\in [K]\\
&B_{i}(s):=A_0 Q_i(s) +\dot Q_i(s)+\sum_{j=1}^K R_{ji}(-\tau_j,s) \qquad i\in [K] \notag \\
&G_{ij}(s,\theta):=\frac{\partial}{\partial s}R_{ij}(s,\theta)+\frac{\partial}{\partial \theta}R_{ji}(s,\theta)^T, \quad i,j\in [K].
\end{align}
Furthermore, suppose
\begin{align}
P&=\tau_KQ_i(0)^T + \tau_KS_i(0) \quad  \text{ for } i\in [K],\\
S_i(s)&=S_i(s)^T,\qquad R_{ij}(s,\theta)=R_{ji}(\theta,s)^T \quad \text{ for } i,j \in  [K],\\
Q_j(s)&=R_{ij}(0,s)\quad  \text{ for } i,j\in [K].
\end{align}
Then the system defined by Equation~\eqref{eqn:delay_eqn} is exponentially stable.
\end{thm}

\begin{IEEEproof} Consider the operator $\mathcal{P}:=\mathcal{P}_{\{P,Q_i,S_i,R_{ij}\}}$.
Since $\{F_1,H_1\}\in \Theta_{n,n,1}$, and  $\left\{M,N\right\} \in \Xi_{d,2n,1}$, by Lemma~\ref{prop:spacing} and Theorem~\ref{thm:positivity_PC}, we have for $x \in Z_{n,K}$
\begin{align}
\ip{x}{\mathcal{P}x}_{Z_{n,K}}-\epsilon \norm{x}^2 =\ip{\hat x}{\left(\mathcal{P}_{M,N}-\mathcal{T}_{F_1,H_1}\right) \hat x}_{L_2^{2n}} =\ip{\hat x}{\mathcal{P}_{M,N} \hat x}_{L_2^{2n}}\ge 0,
\end{align}
where $\hat x \in L_2^{2n}$. This establishes that $\ip{x}{\mathcal{P}x}_{Z_{n,K}} \ge\epsilon \norm{x}^2$ for all $x \in X$. Similarly, examine the operator
\begin{align}
\mathcal{D}\bmat{\bmat{x \\ \phi_1(-\tau_1) \\ \vdots \\ \phi_k(-\tau_K)}\\ \phi_i}(s)=\bmat{\bmat{D_{1} & C_{1} & \cdots &C_{k} \\C_{1}^T &-S_1(-\tau_1) & 0&0 \\
\vdots&0&\ddots&0 \\
C_{k^T} &0&0&-S_k(-\tau_K)}\bmat{x \\ \phi_1(-\tau_1) \\ \vdots \\ \phi_k(-\tau_K)} + \sum_{i=1}^K\int_{-\tau_i}^0 \bmat{B_{i}(s)\\ 0 \\ \vdots \\ 0}\phi_i(s)ds\vspace{3mm} \\
\tau_K B_{i}(s)^T x + \tau_K \dot S_i(s)\phi_i(s)+\sum_{j=1}^K \int_{-\tau_j}^0G_{ij}(s,\theta)\phi_j(\theta) d \theta }.
\end{align}
Since $\{F_2,H_2\}\in \Theta_{n(K+1),n,K}$, and  $\left\{-D,-E\right\} \in \Xi_{d,n(K+2),K}$, we have for $\bmat{x \\ \phi_i} \in Z_{n,K}$
\begin{align}
&\ip{\bmat{\bmat{x \\ \phi_1(-\tau_1) \\ \vdots \\ \phi_k(-\tau_K)}\\ \phi_i}}{\mathcal{D}\bmat{\bmat{x \\ \phi_1(-\tau_1) \\ \vdots \\ \phi_k(-\tau_K)}\\ \phi_i}}_{Z_{n(K+1),n,1}}+\epsilon \norm{\bmat{x \\ \phi_i}}^2 =\ip{z}{\left(\mathcal{P}_{D,E}-\mathcal{T}_{F_2,H_2}\right) z}_{L_2^{n(K+2)}}\\
& =\ip{z}{\mathcal{P}_{D,E} z}_{L_2^{n(K+2)}} \le 0.
\end{align}
This likewise establishes that
\[
\ip{\bmat{\bmat{x \\ \phi_1(-\tau_1) \\ \vdots \\ \phi_k(-\tau_K)}\\ \phi_i}}{\mathcal{D}\bmat{\bmat{x \\ \phi_1(-\tau_1) \\ \vdots \\ \phi_k(-\tau_K)}\\ \phi_i}}_{Z_{n(K+1),n,K}}\le-\epsilon \norm{\bmat{x \\ \phi_i}}^2
\]
for all $\bmat{x\\ \phi} \in X$. By assumption, $P=\tau_KQ_i(0)^T + \tau_KS_i(0)$,  $S_i(s) \in \S^n$, $Q_j(s)=R_{ij}(0,s)$ and $R_{ij}(s,\theta)=R_{ji}(\theta,s)^T$. Hence Theorem~\ref{thm:dual_MD} establishes exponential stability of Equation~\eqref{eqn:delay_eqn}.
\end{IEEEproof}

\section{A Matlab Toolbox Implementation}\vspace{-1mm} \label{sec:toolbox}
To assist with the application of these results, we have created a library of functions for verifying the stability conditions described in this paper. These libraries make use of modified versions of the SOSTOOLS~\cite{prajna_2002} and MULTIPOLY toolboxes coupled with either SeDuMi~\cite{sturm_1999} or Mosek. A complete package can be downloaded from~\cite{mmpeet_web}. Key examples of functions included are:

%
\begin{enumerate}
\item \verb+[M,N]=sosjointpos_mat_ker_ndelay.m+ \vspace{-1mm}
\begin{itemize}
\item Declares a positive piecewise-polynomial multiplier, kernel pair which satisfies $[M,N] \in \Xi_{d,n,K}$.\vspace{-1mm}
\end{itemize}
%
\item \verb+sosmateq.m+ \vspace{-1mm}
\begin{itemize}
\item Declare a matrix-valued equality constraint.\vspace{-1mm}
\end{itemize}

\item \verb+[F,H]=sosspacing_mat_ker_ndelay.m+ \vspace{-1mm}
\begin{itemize}
\item Declare a matrix-valued equality constraint which satisfies $\{F,H\}\in \Theta_{n,n,K}$.\vspace{-1mm}
\end{itemize}
\end{enumerate}
\vspace{2mm}

The functions are implemented within the pvar framework of SOSTOOLS and the user must have some familiarity with this relatively intuitive language to utilize these functions. Note also that the entire toolbox and supporting modified implementations of SOSTOOLS and MULTIPOLY must be added to the path for these functions to execute.

\paragraph{Pseudocode}
To illustrate how these conditions can be efficiently coded using the Matlab toolbox, we give a pseudocode implmentation of the conditions of Theorem~\ref{thm:dualLMI_MD}.
\begin{enumerate}
\item \verb?[M,N]=sosjointpos_mat_ker_ndelay?
\item \verb?[F1,H1]=sosspacing_mat_ker_ndelay?
\item \verb?[D,E]=L(M+F1, N+H1)?
\item \verb?[Q,R]=sosjointpos_mat_ker_ndelay?
\item \verb?[F2,H2]=sosspacing_mat_ker_ndelay?
\item \verb?sosmateq(D+F2+Q)?
\item \verb?sosmateq(E+H2+R)?
\end{enumerate}

Here we use the function $L$ to represent the map $\mathcal{L}_1$. An optimized version of the code is contained in\\
\verb+solver_ndelay_nd_dual_joint.m+.

\section{Numerical Validation}\label{sec:validation}
In this section, we apply the dual stability condition to a battery of numerical examples in order to verify that the proposed stability conditions are not conservative. In each case, the maximum stable value of a specified parameter is given for each degree $d$. In each case $d$ is increased until the maximum parameter value is tight to several decimal places. The computation time is also listed in CPU seconds on an Intel i7-5960X 3.0GHz processor. This time corresponds to the interior-point (IPM) iteration in SeDuMi and does not account for preprocessing, postprocessing, or for the time spent on polynomial manipulations formulating the SDP using SOSTOOLS. Such polynomial manipulations can significantly exceed SDP computation time.

\paragraph{Example A} First, we consider a simple example which is known to be stable for $\tau \le \frac{\pi}{2}$.
\[
\dot x(t)=-x(t-\tau)
\]
\[
\hbox{\begin{tabular}{c|c|c|c|c|c|c}
$d$ & $1$  & $2$ & $3$ & $4$ & $5$ & \text{analytic}\\
\hline
$\tau_{\max}$ & 1.408 & 1.5707 & 1.5707 & 1.5707 & 1.5707& 1.5707\\
\hline
CPU sec & .18 & .21 & .25 & .47 & .73 \\
\end{tabular}}
\]

\paragraph{Example B} Next, we consider a well-studied 2-dimensional, single delay system.
\[
\dot x(t)=\bmat{0 & 1 \\ -2 & .1}x(t)+\bmat{0 & 0\\1 & 0} x(t-\tau)
\]
\[
\hbox{\begin{tabular}{c|c|c|c|c|c}
$d$ & $1$  & $2$ & $3$ & $4$ & \text{limit}\\
\hline
$\tau_{\max}$ & 1.6581 & 1.716 &  1.7178 & 1.7178 &1.7178  \\
$\tau_{\min}$ &  .10019  & .10018 & .10017 & .10017 & .10017\\
\hline
CPU sec & .25 & .344 & .678 & 1.725 & \\
\end{tabular}}
\]

\paragraph{Example C} We consider a scalar, two-delay system.
\[
\dot x(t)=ax(t)+b x(t-1) + c x(t-2)
\]
In this case, we fix $a=-2$, $c=-1$ and search for the maximum $b$, which can be found in, e.g.~\cite{nussbaum_book,gu_2005,egorov_2014} to be 3.
\[
\hbox{\begin{tabular}{c|c|c|c|c|c}
$d$ & $1$  & $2$ & $3$ & $4$ & \text{analytic}\\
\hline
$b_{\max}$ & .7071 & 2.5895 & 2.9981 & 2.9982 & 3\\
\hline
CPU sec & .3 & .976 & 2.77 & 12.96 & \\
\end{tabular}}
\]

\paragraph{Example D} We consider a 2-dimensional, two-delay system where $\tau_1=\tau_2/2$ and search for the maximum stable $\tau_2$.
\[
\dot x(t)=\bmat{0 & 1\\ -1 & .1}x(t)+\bmat{0 & 0\\-1 & 0} x(t-\tau/2) + \bmat{0 & 0\\1 & 0} x(t-\tau)
\]
\[
\hbox{\begin{tabular}{c|c|c|c|c|c}\label{tab:taumax}
$d$ & $1$  & $2$ & $3$ & $4$ & \text{limit}\\
\hline
$\tau_{\max}$ & 1.33 & 1.371 & 1.3717 & 1.3718 & 1.372\\
\hline
CPU sec & 2.13 & 6.29 & 24.45 & 79.0 & \\
\end{tabular}}
\]

\section{Conclusion}\label{sec:conclusion}
In conclusion, we have proposed a new form of dual Lyapunov stability condition which allows convexification of the controller synthesis problem for delayed and other infinite-dimensional systems. This dual principle requires a Lyapunov operator which is positive, invertible, self-adjoint and preserves the structure of the state-space. We have proposed such a class of operators and used them to create stability conditions which can be expressed as positivity and negativity of quadratic Lyapunov functions. These dual stability conditions have a tridiagonal structure which is distinct from standard Lyapunov-Krasovskii forms and may be exploited to increase performance when studying systems with large numbers of delays. The dual stability condition is presented in a format which can be adapted to many existing computational methods for Lyapunov stability analysis. We have applied the Sum-of-Squares approach to enforce positivity of the quadratic forms and tested the stability condition in both the single and multiple-delay cases. Numerical testing on several examples indicates the method is not conservative. The contribution of the present paper is not in the efficiency of the stability test, however, as these are likely less efficient when compared to, e.g., previous SOS results due to the highly structured nature of the operators used. Rather the contribution is in the convexification of the synthesis problem which opens the door for dynamic output-feedback $H_\infty$ synthesis for infinite-dimensional systems.

\appendices
\section{Table of Notation}
For convenience, we summarize a selected subset of the notation used in this paper.

\vspace{3mm}
\noindent \textbf{Spaces:} $Z_{m,n,K}:=\{\R^m \times L_2^{n}[-\tau_1,0]\times \cdots \times L_2^n[-\tau_K,0]\}$
with $Z_{n,K}:=Z_{n,n,K}$ and
\[
\ip{\bmat{y\\ \psi_i}}{\bmat{x\\ \phi_i}}_{Z_{m,n,K}}=\tau_K y^T x + \sum_{i=1}^K \int_{-\tau_i}^0 \psi_i(s)^T\phi_i(s)ds.
\]

\vspace{3mm}
\noindent \textbf{Subsets:}

\[
\Xi_{d,n,K}:=\{\{M,N\}\,:\,  \substack{ M=M_1+M_2,\, N=N_1+N_2,\, \text{ where $\{M_1,N_1\}$ and $\{M_2,N_2\}$ satisfy the}\\
\text{conditions of Thm.~\ref{thm:positivity_PC} with $g_i=1$ and $g_i=-(s+\tau_i)(s+\tau_{i-1})$, respectively.}}  \}
\]

\[
\Theta_{m,n,K}:=\{\{F,H\}\,:\, F,H \text{ satisfy the conditions of Thm.~\ref{thm:mixed_pos}.}\}
\]

\vspace{3mm}
\noindent \textbf{Operators:}
\begin{align}
&\left(\mathcal{P}_{\{P,Q_i,S_i,R_{ij}\}} \bmat{x\\ \phi_i}\right)(s) :=\bmat{  P x +  \sum_{i=1}^K \int_{-\tau_i}^0 Q_i(s)\phi_i(s) d s \\
\tau_K Q_i(s)^T x + \tau_K S_i(s)\phi_i(s) + \sum_{j=1}^K \int_{-\tau_j}^0 R_{ij}(s,\theta)\phi_j(\theta)\, d \theta }.
\end{align}

\begin{equation}
\left(\mathcal{P}_{M,N}x\right)(s):= M(s)x(s) + \int_{-\tau_K}^0 N(s,\theta)x(\theta)d \theta.\vspace{-1mm}
\end{equation}

Given $\{F,H\} \in \Theta_{m,n,K}$, we define the operator $\mathcal{T}_{F,H}: L_2^{m,n}\rightarrow L_2^{m,n}$ as
\[
\left(\mathcal{T}_{F,H}z\right)(s):=F(s)z(s)+\int_{-\tau_K}^0H(s,\theta)z(\theta)\, d\theta.
\]

\vspace{3mm}
\noindent \textbf{Linear Transformations:} We say
\[
\{M,N\}:=\mathcal{L}_1(P,Q_i,S_i,R_{ij})
\]
if $a_i=\frac{\tau_i-\tau_{i-1}}{\tau_i}$ and
\begin{align}
M(s)&=\begin{cases}
\bmat{P & \frac{\tau_K}{a_i} Q_i(\frac{s+\tau_{i-1}}{a_i})\\
\frac{\tau_K}{a_i}  Q_i(\frac{s+\tau_{i-1}}{a_i})^T & \frac{\tau_K}{a_i} S_i(\frac{s+\tau_{i-1}}{a_i}) } & s \in [-\tau_i, -\tau_{i-1}]\\
\end{cases}\notag \\
N(s,\theta)&=\begin{cases}
R_{ij}(\frac{s+\tau_{i-1}}{a_i},\frac{\theta+\tau_{j-1}}{a_j}) & s \in [-\tau_i, -\tau_{i-1}],\, \theta \in [-\tau_j, -\tau_{j-1}]\\
\end{cases}
\end{align}

\section*{Acknowledgment}
This work was supported by the National Science Foundation under grants No. 1100376 and 1301851.

\bibliographystyle{IEEEtran}
\bibliography{peet_bib,delay,NSF_CAREER_bib2011,LMIs}

\begin{IEEEbiography}{Matthew M. Peet}
received the B.S. degree in physics and in aerospace engineering from the University of Texas, Austin, TX, USA, in 1999 and the M.S. and Ph.D. degrees in aeronautics and astronautics from Stanford University, Stanford, CA, in 2001 and 2006, respectively. He was a Postdoctoral Fellow at the National Institute for Research in Computer Science and Control (INRIA), Paris, France, from 2006 to 2008, where he worked in the SISYPHE and BANG groups. He was an Assistant Professor of Aerospace Engineering in the Mechanical, Materials, and Aerospace Engineering Department, Illinois Institute of Technology, Chicago, IL, USA, from 2008 to 2012. Currently, he is an Assistant Professor of Aerospace Engineering, School for the Engineering of Matter, Transport, and Energy, Arizona State University, Tempe, AZ, USA, and Director of the Cybernetic Systems and Controls Laboratory. Dr. Peet received a National Science Foundation CAREER award in 2011.
\end{IEEEbiography}

\end{document}